\documentclass[11pt]{article}
\usepackage{epsf}
\usepackage{color}
\usepackage{epsfig}

\usepackage{graphicx}
\usepackage{amsmath}
\usepackage{latexsym}
\usepackage{amssymb}

\usepackage{a4}

\newtheorem{theorem}{Theorem}

\newcommand{\e}{{\rm e}}

\renewcommand{\r}{\mathbb{R}}

\newcommand{\cO}{{\mathcal O}}

\newcommand{\cM}{\r^D}

\newcommand{\cL}{{\mathcal L}}

\newcommand{\cinf}{C^{\infty}(\cM,\r)}
\newcommand{\psib}{\chi}

\newcommand{\psis}{{\mathcal S}}

\begin{document}

\title{Splitting and composition methods in the numerical integration of differential equations}

\author{Sergio Blanes$^{1}$\thanks{Email: \texttt{serblaza@imm.upv.es}}
   \and
%    author two information
  Fernando Casas$^{2}$\thanks{Email: \texttt{Fernando.Casas@uji.es}}
   \and
Ander Murua$^{3}$\thanks{Email: \texttt{Ander.Murua@ehu.es}}
   }

\date{}
\maketitle

\begin{abstract}

We provide a comprehensive survey of splitting and composition methods for the numerical
integration of ordinary differential equations (ODEs). Splitting methods constitute an appropriate choice
when the vector field associated with the ODE can be decomposed into several pieces and each of them
is integrable. This class of integrators are explicit, simple to implement and preserve
structural properties of the system. In consequence, they are specially useful in geometric numerical
integration. In addition, the numerical solution obtained by splitting schemes can be seen as the
exact solution to a perturbed system of ODEs possessing the same geometric properties as the
original system. This backward error interpretation has direct implications for the qualitative
behavior of the numerical solution as well as for the error propagation along time. Closely
connected with splitting integrators are composition methods. We analyze the order conditions
required by a method to achieve a given order and summarize the different families of schemes one
can find in the literature. Finally, we illustrate the main features of splitting and composition methods
on several numerical examples arising from applications.

\vspace*{0.5cm}

\begin{description}
 \item $^1$Instituto de Matem\'atica Multidisciplinar,
  Universidad Polit\'ecnica de Valencia, E-46022  Valencia, Spain.
 \item $^2$Departament de Matem\`atiques, Universitat Jaume I,
  E-12071 Castell\'on, Spain.
 \item $^3$Konputazio Zientziak eta A.A. saila, Informatika
Fakultatea, EHU/UPV, Donostia/San Sebasti\'an, Spain.
\end{description}

\end{abstract}

\section{Introduction by examples}

The basic idea of splitting methods for the time integration of ordinary differential equations (ODEs)
can be formulated as follows. Given the initial value problem
\begin{equation}   \label{eq.1.1}
   x' = f(x), \qquad x_0 = x(0) \in \mathbb{R}^D
\end{equation}
with $f: \mathbb{R}^D \longrightarrow  \mathbb{R}^D$ and solution $\varphi_t(x_0)$,
let us suppose that  $f$  can be
expressed as $f = \sum_{i=1}^m f^{[i]}$ for certain functions
$f^{[i]}: \mathbb{R}^D \longrightarrow  \mathbb{R}^D$, in such a way that the equations
\begin{equation}   \label{eq.1.2}
   x' = f^{[i]}(x), \qquad x_0 = x(0) \in \mathbb{R}^D, \qquad i=1, \ldots, m
\end{equation}
can be integrated exactly, with solutions $x(h) = \varphi_h^{[i]}(x_0)$ at $t = h$, the
time step. Then, by combining these solutions as
\begin{equation}   \label{eq.1.3}
   \psib_h = \varphi_h^{[m]} \circ \cdots \circ \varphi_h^{[2]} \circ \varphi_h^{[1]}
\end{equation}
and expanding $\psib$ into Taylor series, one finds that $\psib_h(x_0) =
\varphi_h(x_0) + \mathcal{O}(h^2)$, so that $\psib_h$ provides a
first-order approximation to the exact solution. As
we will see, it is possible to get higher order approximations by
introducing more maps with additional coefficients,
$\varphi_{a_{ij} h}^{[i]}$,  in the previous composition
(\ref{eq.1.3}).

One thus may say that splitting methods involve three steps: (i) choosing the set
of functions $f^{[i]}$ such that $f = \sum_i f^{[i]}$; (ii)
solving either exactly or approximately each equation $x' =
f^{[i]}(x)$; and (iii) combining these solutions to construct an
approximation for  (\ref{eq.1.1}). One obvious requirement is that
the equations $x' = f^{[i]}(x)$ should be simpler to integrate
than the original system.

Some of the advantages that splitting methods possess can be summarized as follows:
\begin{itemize}
\item They are usually simple to implement.
\item They are, in general, explicit.
\item Their storage requirements are quite modest. The algorithms are sequential
  and the solutions at intermediate stages are stored in the
  solution vectors. This property can be of great interest when they are applied to
  partial differential equations (PDEs) previously semidiscretized.
\item There exist in the literature a large number of specific methods
   tailored for different structures.
\item They preserve structural properties of the exact solution,
   thus conferring to the numerical scheme a qualitative superiority
   with respect to other standard integrators, especially when long
   time intervals are considered. Examples of these structural features
   are symplecticity, volume preservation, time-symmetry and
   conservation of first integrals. In this sense, splitting methods constitute an
   important class of \emph{geometric numerical integrators}.
\end{itemize}

Let us give more details on this last item.
Traditionally, the goal of numerical integration of ODEs consists
in computing the solution to the initial value problem
(\ref{eq.1.1}) at time $t_N = N h$ with a global error $\|x_N -
x(t_N)\|$ smaller than a prescribed tolerance and as efficiently
as possible. To do that one chooses the class of method (one-step,
multistep, extrapolation, etc.), the order (fixed or adaptive) and the time step
(constant or variable). In contrast, with a geometric numerical
integrator one typically fix a (not necessarily small) time step
and compute solutions for very long times for several initial
conditions, in order to get an approximate phase portrait of the
system. It turns out that, although the global error of each
trajectory may be large, the phase portrait is in some sense close
to that of the original system.

The aim of geometric numerical integration is thus to reproduce the
qualitative features of the solution of the differential equation
which is being discretised, in particular its geometric
properties \cite{budd99gin,hairer06gni}. The motivation for developing such
structure-preserving algorithms arises independently in areas of
research as diverse as celestial mechanics, molecular dynamics,
control theory, particle accelerators physics, and numerical
analysis
\cite{hairer06gni,iserles00lgm,mclachlan02sm,mclachlan06gif,leimkuhler04shd}.
Although diverse, the systems appearing in these areas have
one important common feature. They all preserve some
underlying geometric structure which influences
the qualitative nature of the phenomena they produce.
In the field of geometric numerical integration these properties are built into
the numerical method, which gives the method an improved
qualitative behaviour, but also allows for a significantly more
accurate long-time integration than with general-purpose methods.
In addition to the construction of new numerical algorithms,
an important aspect of geometric integration is the explanation
of the relationship between preservation of the geometric
properties of the scheme and the observed favorable
error propagation in long-time integration.

Before proceeding further, let us introduce at this point some splitting methods
and illustrate them on simple examples.

\paragraph{Example 1: Symplectic Euler and leapfrog schemes.}

Suppose we have a Hamiltonian system of the form $H(q,p) = T(p) + V(q)$, where $q \in \mathbb{R}^d$
are the canonical coordinates, $p \in \mathbb{R}^d$ are the conjugate momenta, $T$ represents the kinetic
energy and $V$ is the potential energy. Then the equations of motion read \cite{goldstein80cme}
\begin{equation}  \label{eq.1.4}
   q' = T_p(p), \qquad p' = -V_q(q),
\end{equation}
where $T_p$ and $V_q$ denote the vectors of partial derivatives.  Equations (\ref{eq.1.4})
can be formulated as (\ref{eq.1.1}) with $x=(q,p)^T$, $f(x) = (T_p,-V_q)^T = J \nabla H(x)$ and
$D=2d$. Here $J$ denotes the $2d \times 2d$ canonical symplectic matrix
\[
   J = \left( \begin{array}{rr}
               0  &  I_d  \\
              -I_d  &  0
             \end{array}  \right)
\]
and $I_d$ stands for the $d$-dimensional identity matrix. In this case the exact flow $\varphi_t$ is
symplectic \cite{arnold89mmo}.
The simple Euler method applied to this system
provides the following first order approximation for a time step $h$:
\begin{equation}   \label{Euler}
  \begin{array}{lcl}
   q_{n+1} & = & q_n + h T_p(p_n)   \\
   p_{n+1} & = & p_n - h V_q(q_n).
   \end{array}
\end{equation}
On the other hand, if we consider $H$ as the sum of two
Hamiltonians, the first one depending only on $p$ and the second only on $q$, the corresponding
Hamilton equations
\begin{equation}   \label{eq.1.5}
    \begin{array}{lcl}
       q' & = & T_p(p) \\
       p' & = &  0
    \end{array}
    \qquad  \mbox{ and } \qquad
    \begin{array}{lcl}
       q' & = & 0 \\
       p' & = &  -V_q(q)
    \end{array}
\end{equation}
with initial condition $(q_0,p_0)$ can be readily solved to yield
\begin{equation}   \label{eq.1.6}
  \varphi_t^{[T]}:  \; \begin{array}{lcl}
      q(t) & = & q_0 + t \, T_p(p_0)  \\
      p(t) & = & p_0
    \end{array}
    \qquad  \mbox{ and } \qquad
   \varphi_t^{[V]}:  \; \begin{array}{lcl}
      q(t) & = & q_0  \\
      p(t) & = & p_0 - t \, V_q(q_0),
    \end{array}
\end{equation}
respectively. Composing the time $t = h$ flow $\varphi_h^{[V]}$ (from initial condition $(q_n,p_n)$) followed
by $\varphi_h^{[T]}$, gives the scheme
\begin{equation}   \label{Euler-sympl}
  \chi_h \equiv \varphi_h^{[T]} \circ \varphi_h^{[V]}:   \;  \begin{array}{lcl}
      p_{n+1} & = & p_n - h \, V_q(q_{n}) \\
      q_{n+1} & = & q_n + h \, T_p(p_{n+1}).
    \end{array}
\end{equation}
Since it is a composition of the flows of two Hamiltonian systems and in addition
the composition of symplectic
maps is again symplectic, $\chi_h$ is a symplectic integrator, which can be called
appropriately the \emph{symplectic Euler method}. It is of course also possible to compose the maps in the
opposite order, $\varphi_h^{[V]} \circ  \varphi_h^{[T]}$, thus obtaining another first order symplectic Euler scheme:
\begin{equation}   \label{Euler-sympl2}
  \chi_h^* \equiv \varphi_h^{[V]} \circ \varphi_h^{[T]}:   \;  \begin{array}{lcl}
      q_{n+1} & = & q_n + h \, T_p(p_n)  \\
      p_{n+1} & = & p_n - h \, V_q(q_{n+1}).
    \end{array}
\end{equation}
One says that (\ref{Euler-sympl2}) is the \emph{adjoint} of $\chi_h$. Yet another possibility consists in using a `symmetric' version
\begin{equation}  \label{leapfrog}
  \mathcal{S}_h^{[2]} \equiv \varphi_{h/2}^{[V]} \circ \varphi_h^{[T]} \circ  \varphi_{h/2}^{[V]},
\end{equation}
which is known as the Strang splitting \cite{strang68otc}, the leapfrog or the St\"ormer--Verlet method
\cite{verlet67ceo}, depending on the context where it
is used. Observe that $\mathcal{S}_h^{[2]} = \chi_{h/2} \circ \chi_{h/2}^*$ and it is also symplectic and
second order.

\paragraph{Example 2: Harmonic oscillator.}

Let us consider now the Hamiltonian function $H(q,p) = \frac{1}{2}(p^2 + q^2)$, where now
$q,p \in \mathbb{R}$. Then the
corresponding equations (\ref{eq.1.4}) are linear  and can be written as
\begin{equation}\label{harmonic2}
    x' \equiv \left(\begin{array}{c}
        q^\prime \\
        p^\prime  \end{array}\right) =
    \Big[   \underbrace{\left(\begin{array}{cc}
        0 & 1 \\
        0 & 0  \end{array}\right)}_{A} +
         \underbrace{\left(\begin{array}{cc}
        0 & 0 \\
       -1 & 0  \end{array}\right)}_{B}
 \Big]
     \left(\begin{array}{c}
        q \\
        p  \end{array}\right) = (A + B) \, x.
\end{equation}
This system has periodic solutions for which the energy $H$ is conserved. In addition, it is area
preserving and time reversible.
The numerical solution obtained by the Euler scheme (\ref{Euler}) reads
\begin{equation}   \label{eq.1.7}
   \left(\begin{array}{c}
        q_{n+1} \\
        p_{n+1}  \end{array}\right) =
        \left(  \begin{array}{cc}
           1  &  h  \\
           -h & 1
             \end{array} \right) \,
                \left(\begin{array}{c}
        q_{n} \\
        p_{n}  \end{array}\right),
\end{equation}
whereas the symplectic Euler method (\ref{Euler-sympl2}) leads to
\begin{equation}  \label{eq.1.8}
   \left(\begin{array}{c}
        q_{n+1} \\
        p_{n+1}  \end{array}\right) =
        \left(  \begin{array}{cc}
           1  &  h  \\
           -h & 1-h^2
             \end{array} \right) \,
                \left(\begin{array}{c}
        q_{n} \\
        p_{n}  \end{array}\right) = \e^{h B} \e^{h A}
              \left(\begin{array}{c}
        q_{n} \\
        p_{n}  \end{array}\right).
\end{equation}
Both render first order approximations to the exact solution, which can be expressed as
$x(t) = \e^{h(A+B)} x_0$,
but there are important differences between them. First, the map (\ref{eq.1.8}) is area preserving (because it
is symplectic), in contrast with (\ref{eq.1.7}). Second, the approximation obtained by the symplectic Euler scheme
verifies
\[
   \frac{1}{2} (p_{n+1}^2 + h p_{n+1} q_{n+1} + q_{n+1}^2) = \frac{1}{2} (p_{n}^2 + h p_{n} q_{n} + q_{n}^2).
\]
Third, it can be shown that (\ref{eq.1.8}) \emph{is} the exact solution at $t = h$ of
the \emph{perturbed} Hamiltonian system
\begin{eqnarray}  \label{eq.1.9}
  \tilde{H}(q,p,h) &=& \frac{2\arcsin(h/2)}{h\sqrt{4-h^2}} (p^2 + h \, p \, q + q^2)\\
 & =  & \frac{1}{2} (p^2 + q^2) + h\, \left(\frac12 \, p \, q + \frac{1}{12} h (p^2 + q^2) + \cdots \right). \nonumber
\end{eqnarray}
In other words, the numerical approximation (\ref{eq.1.8}), which is only
of first order for the exact trajectories of the Hamiltonian $H(q,p) = \frac{1}{2}(p^2 + q^2)$, is in fact the
exact solution of the perturbed Hamiltonian (\ref{eq.1.9}).

\begin{figure}[h!]
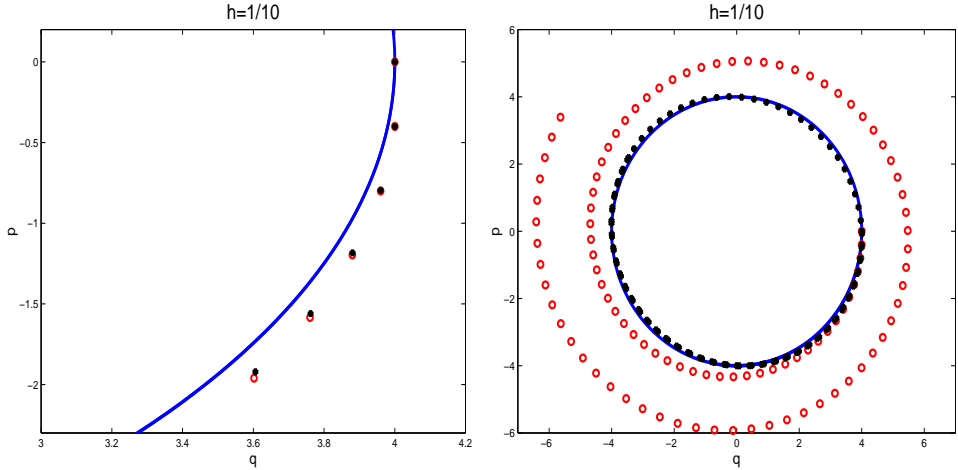

\begin{center}
\makebox{\psfig{figure=Harm5pasos.eps,height=6.3cm,width=6.3cm}}
\makebox{\psfig{figure=Harm100pasos.eps,height=6.3cm,width=6.3cm}}
\end{center}
\caption{{Numerical integration of the harmonic oscillator using
the Euler method (white circles) and the symplectic Euler method
(black circles) with initial condition $(q_0,p_0)=(4,0)$ and  time
step $h=\frac{1}{10}$. The left panel shows the results for the
first 5 steps, whereas the right panel shows the results for the
first 100 steps. The exact solution corresponds to the solid
line.}}
 \label{fig1}
\end{figure}

How these features manifest in actual simulations? To illustrate
this point we take initial conditions $(q_0,p_0)=(4,0)$ and integrate
with a time step $h=0.1$. Figure~\ref{fig1} shows the first five
numerical approximations obtained by the Euler method
(\ref{eq.1.7}) and the symplectic Euler scheme (\ref{eq.1.8}) in the left panel, and
the results for the first 100 steps in the right panel. It is clear that for one
time step there are not significant differences between the standard
Euler and the symplectic Euler methods, but the picture is completely
different for longer integrations, where the superiority of the
splitting symplectic method is evident. Note that  the numerical solution it provides
evolves on a slightly perturbed ellipse.

\paragraph{Example 3: The 2-body problem (Kepler problem).}

The motion of two bodies attracting each other through the gravitational law can be described by
\begin{equation}   \label{eq.1.10}
   q_i^{\prime\prime} = - \frac{q_i}{(q_1^2 + q_2^2)^{3/2}}, \qquad i=1,2
\end{equation}
in conveniently normalized coordinates in the plane of motion. This system has a number
of characteristic geometric properties. First, equations (\ref{eq.1.10})
can be derived from the Hamiltonian function
\begin{equation}  \label{eq.1.11}
  H(q,p) = T(p) + V(q) = \frac{1}{2} (p_1^2 + p_2^2) - \frac{1}{r}, \qquad r = \sqrt{q_1^2 + q_2^2}.
\end{equation}
Second, it is invariant under continuous translations in time and rotations in space, and thus both the
Hamiltonian and the angular momentum $L = q_1 p_2 - q_2 p_1$ are preserved. In addition, the
so-called Laplace--Runge--Lenz vector is also constant along their solutions, due to
the fact that the symmetry group of this problem is the group of four-dimensional real proper rotations
$\mathrm{SO}(4)$ \cite{goldstein80cme}.

For the numerical integration of this problem we choose as
initial value
\begin{equation}\label{eq.1.12}
  q_1(0) = 1- e, \quad q_2(0) = 0, \quad p_1(0) = 0, \quad p_2(0) = \sqrt{\frac{1+e}{1-e}},
\end{equation}
where $0 \le e < 1$ represents the eccentricity of the orbit. In
this case the total energy is $H = H_0 = -1/2$, the period of the
solution is $2 \pi$, the initial condition corresponds to the
pericenter and the major semiaxis of the ellipse is 1.

\begin{figure}[h!]
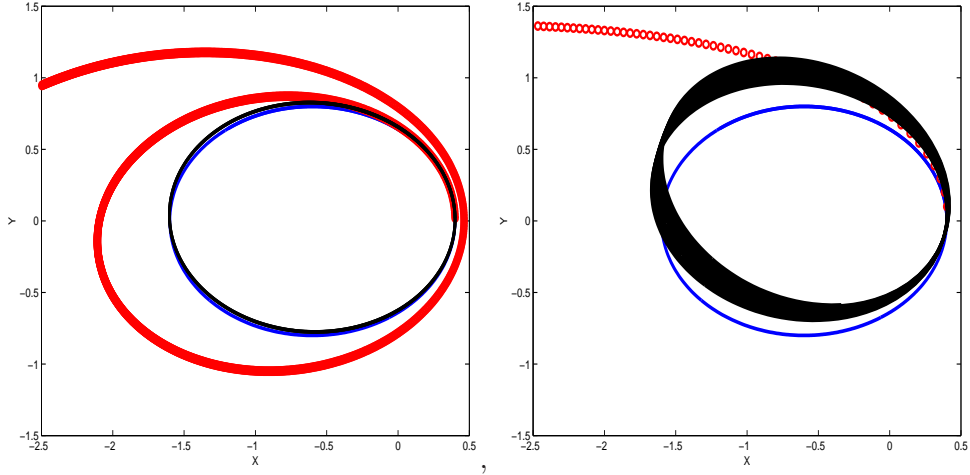

\begin{center}
\makebox{\psfig{figure=Kep_h01Per3.eps,height=6.3cm,width=6.3cm}},
\makebox{\psfig{figure=Kep_h05Per15.eps,height=6.3cm,width=6.3cm}}
\end{center}
\caption{{Numerical integration of the 2-body problem using the
Euler method (white circles) and the symplectic Euper scheme
(black circles) for the initial conditions (\ref{eq.1.12}) with
eccentricity $e=0.6$. The left panel shows the results for
$h=\frac{1}{100}$ and the first 3 periods and the right panel
shows the results for $h=\frac{1}{20}$ and the first 15 periods.}}
 \label{fig2}
\end{figure}

Figure~\ref{fig2} shows some numerical solutions obtained with
schemes (\ref{Euler}) and (\ref{Euler-sympl})  for the initial
conditions (\ref{eq.1.12}) with eccentricity $e=0.6$. The left panel shows the
results for the integration of  3 periods with time step
$h=\frac{1}{100}$.  As in the previous example, the explicit Euler
method provides an approximate orbit that spirals outwards,
whereas the symplectic Euler scheme merely distorts the ellipse, but also
exhibits a precession effect.
To better illustrate this effect, we
repeat the experiment for a longer interval (15 periods) and a
larger time step ($h=\frac{1}{20}$) in the right panel. The explanation of these phenomena can
be formulated as follows. On the one hand, the symplectic Euler method
exactly conserves the angular momentum. On the other hand, the numerical solution it provides
can be seen as the exact solution of a
slightly perturbed Kepler problem, and thus $\mathrm{SO}(4)$ is no longer the symmetry group
of the problem, so that the Laplace--Runge--Lenz vector is not preserved and the trajectories are not closed
anymore.

\begin{figure}[h!]
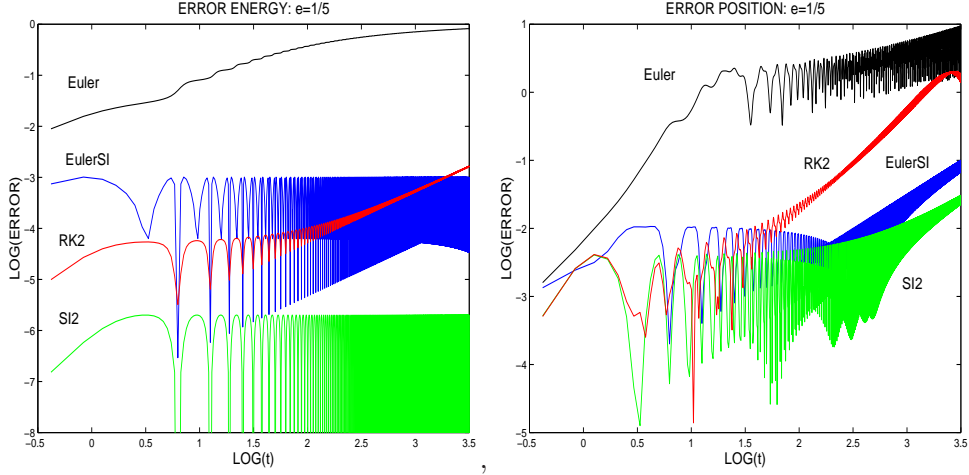

\begin{center}
\makebox{\psfig{figure=KepErrorEn2.eps,height=6.3cm,width=6.3cm}},
\makebox{\psfig{figure=KepErrorPos2.eps,height=6.3cm,width=6.3cm}}
\end{center}
\caption{{Error growth in energy and position for the Kepler
problem with $e=1/5$ and 500 periods achieved by
the first order symplectic Euler (EulerSI) and the second order St\"ormer--Verlet
integrator (SI2). For comparison, we also include
the first order Euler and the second order Heun (RK2) methods.
The time step is adjusted in such a way
that all methods use 1500 force evaluations per period.}}
 \label{fig3}
\end{figure}

Next we check how the error in the preservation
of energy and the global error in position propagates with
time. For comparison, we also include the results obtained with a Runge--Kutta
method of order 2 (Heun's method) and the leapfrog scheme
(\ref{leapfrog}). We now consider $e=1/5$ and integrate for 500
periods. The step size is $h=\frac{2\pi}{1500}$ in all cases, except for the Heun
method which uses  $h=\frac{2\pi}{750}$ instead. In this way all methods require
the same number of force evaluations (Heun's method computes twice the force per step).
The corresponding results are shown in Figure \ref{fig3} in a log-log scale.
Notice that   the average error in energy does not grow for the split
symplectic methods and the error in positions  grows only linearly with time, in contrast with
Euler and Heun schemes. The  St\"ormer--Verlet  integrator provides more accurate results than
the Heun method with the same computational cost.

\paragraph{A collection of (additional) examples.}

Splitting methods constitute an important tool in different areas of science.
In addition to Hamiltonian systems,
they can be successfully applied in the numerical study of Poisson systems,
systems possessing integrals
of the motion (such as energy and angular momentum) and systems
with (continuous, discrete and time-reversal)
symmetries. As a matter of fact, splitting methods have been designed
(often independently) and extensively used
in fields as distant as molecular dynamics, simulation of storage rings
in particle accelerators, celestial
mechanics and quantum physics simulations. To see why this is so, next we collect
a number of differential equations which appear in different contexts ranging from
Celestial Mechanics to electromagnetism and Quantum
Mechanics. These examples also try to illustrate the fact that very often one
particular equation can be split into different ways and the most
appropriate methods may depend on the split chosen.
%This
%will also motivate the search of methods for special problems.

We (arbitrarily) classify our examples into three different categories.

\begin{enumerate}

\item Hamiltonian systems.
\begin{description}
\item[(a)] Generalized harmonic oscillator ($M,N\in\mathbb{R}^{d\times
d}$):
 \begin{equation}  \label{harmonicGeneral}
   H=\frac12 p^TMp + \frac12 q^TNq.
 \end{equation}

\item[(b)] H\'{e}non--Heiles Hamiltonian \cite{henon64tao}:
 \begin{equation}  \label{HenonHeiles}
  H=\frac12 (p_{1}^{2}+p_{2}^{2})+\frac12
  (q_1^{2}+q_2^{2})+q_1^{2}q_2-\frac13 q_2^{3}.
 \end{equation}

\item[(c)] Perturbed Kepler problem. It models the dynamics of a satellite moving into the
  gravitational field produced by a slightly oblate planet:
 \begin{equation}   \label{eq.14}
   H = \frac{1}{2} (p_1^2 + p_2^2) - \frac{1}{r} -
       \frac{\varepsilon}{2 r^3} \left( 1 - \frac{3 q_1^2}{r^2}
       \right)
 \end{equation}
where $\varepsilon$ is typically a small parameter. When $\varepsilon = 0$, the 2-body problem
(\ref{eq.1.11}) is recovered.

\item[(d)] The gravitational $N$-body problem ($q_i,p_i\in\mathbb{R}^3, \
i=1,\ldots,N$):
 \begin{equation}   \label{eq.14N}
   H = \frac{1}{2} \sum_{i=1}^N  \frac{1}{2m_i}p_i^Tp_i -
    G \sum_{i=2}^N \sum_{j=1}^{i-1}
    \frac{m_im_j}{\|q_i-q_j\|}.
 \end{equation}

\item[(e)] The motion of a charged particle in a constant
 magnetic field perturbed by $k$ electrostatic plane waves \cite{candy91asi}:
 \begin{equation}\label{wave1}
  H(q,p,t) = \frac12p^{2}+\frac12q^{2}+
          \varepsilon \sum_{i=1}^k\cos (q-\omega_it).
 \end{equation}

\end{description}

\item More general dynamical systems.

\begin{description}
  \item[(a)] The Volterra--Lotka problem \cite{hairer06gni},
\begin{equation} \label{volterra-lotka}
  \frac{d}{dt} \left[
\begin{array}{c}
  u  \\
  v
\end{array} \right] =
 \left[
\begin{array}{rr}
   -2 & 0  \\
   0   &  1
\end{array} \right]
 \left[
\begin{array}{c}
  u  \\
  v
\end{array} \right] +
 \left[
\begin{array}{r}
   u v \\
   -u v
\end{array} \right] =
 \left[
\begin{array}{c}
  u(v-2)  \\
  0
\end{array} \right]+
 \left[
\begin{array}{c}
  0  \\
  v(1-u)
\end{array} \right],
\end{equation}
with first integral $I(u,v) = \log u - u + 2 \log v - v$.

  \item[(b)] The Lorenz system \cite{lorenz63dnf,guckenheimer83nod} (split into linear and non-linear parts):
   \begin{equation}
  \frac{d}{dt} \left[
\begin{array}{c}
  x  \\
  y  \\
  z
\end{array} \right] =
 \left[
\begin{array}{rrr}
  -\sigma & \sigma &  0  \\
       r  &   -1   &  0   \\
       0  &    0   &  b
\end{array} \right]
 \left[
\begin{array}{c}
  x  \\
  y  \\
  z
\end{array} \right] +
 \left[
\begin{array}{rrr}
   0  &  0 &  0  \\
   0  &  0 & -x   \\
   0  &  x &  0
\end{array} \right]
 \left[
\begin{array}{c}
  x  \\
  y  \\
  z
\end{array} \right] .
\end{equation}
Here $\sigma, r, b > 0$ are constant. The values considered in \cite{lorenz63dnf} were
$\sigma = 10$, $r=28$ and $b=8/3$.

  \item[(c)] The ABC-flow ($f =f_A + f_B + f_C$, but other splits are
also possible \cite{mclachlan02sm}):
\begin{equation}
 \frac{d}{dt} (x,y,z)
   = A(0,\sin x, \cos x) + B(\cos y , 0 , \sin y)
   + C(\sin z, \cos z , 0).
\end{equation}

\end{description}

\item Evolutionary PDEs.

Although we are mainly concerned here with splitting methods applied
to ODEs, it turns out that they can also be used in the numerical integration of certain partial
differential equations. Specifically, a number of PDEs relevant in the applications, after an
appropriate space discretization, lead to a system of ODEs which can be subsequently solved
numerically by splitting methods.  Among these equations, the following are worth to be mentioned.
\begin{description}
\item[(a)] The Schr\"odinger equation ($\hbar = 1$):
 \begin{equation}\label{Schrodinger}
  i\frac{\partial}{\partial t} \Psi(x,t) =
  \left( -\frac{1}{2m} \nabla^2 + V(x) \right)  \Psi(x,t).
 \end{equation}

\item[(b)] The Gross--Pitaevskii equation \cite{pitaevskii03bec}:
\begin{equation}\label{Gross-Pitaevskii}
  i\frac{\partial}{\partial t} \Psi(x,t) =
  \left( -\frac{1}{2m} \nabla^2  + V(x)  +
    \alpha |\Psi(x,t)|^2 \right) \Psi(x,t)
\end{equation}

\item[(c)] The Maxwell equations
 \begin{equation}\label{Maxwell}
  \frac{\partial}{\partial t} {\bf B} = -\frac{1}{\mu}
 {\bf \nabla \times E}, \qquad
 \frac{\partial}{\partial t} {\bf E} = \frac{1}{\varepsilon}
 {\bf \nabla \times B},
 \end{equation}
 where ${\bf E}(x,t)$, ${\bf B}(x,t)$ are the electric and magnetic
 field vectors, $\mu(x)$ is the  the permeability and
$\varepsilon(x)$ is the permittivity.

\end{description}

\end{enumerate}

\paragraph{}

As we stated before, splitting methods form a subclass of geometric numerical integrators for various
types of ODEs. The reason is easy to grasp from the examples analyzed before.
Suppose that the flow of the original differential equation
(\ref{eq.1.1}) forms a particular group of diffeomorphisms
(in the case of Hamiltonian system, the group of symplectic maps). If $f$ is
conveniently split as $f = \sum_i f^{[i]}$ (step (i) in the construction process of a splitting scheme) and the flows corresponding to each $f^{[i]}$ also belong to the same group
of diffeomorphisms  in such a way that they can be explicitly obtained (step (ii)), then, by
composing these flows (step (iii)) we get an approximation in the group, thus inheriting geometric properties
of the exact solution. These considerations also hold (with some modifications) when the exact flow forms
a semigroup or a symmetric space.

With respect to steps (i) and (ii) before, several  comments are
in order. First, whereas for certain classes of ODEs, the
splitting can be constructed systematically for any $f$, in other
cases no general procedure is known, and thus one has to proceed
on a case by case basis. Second, sometimes a standard splitting is
possible for a certain $f$, but there exist other possible choices
leading to more efficient schemes (we will see some examples in
section \ref{nexamples}). Third, whereas the original system possesses several
geometric properties which are interesting to preserve  by the
numerical scheme, different splittings preserve different
properties and it is not always possible to find one splitting
preserving all of them. These aspects have been analyzed in detail
in \cite{mclachlan02sm}, where a classification of ODEs and general
guidelines to find suitable splittings in each case is provided.
Here, by contrast, we will concentrate on the third step of any
splitting method: given a particular splitting, we will show how
to combine the flows of the pieces $f^{[i]}$ to get efficient
higher order approximations. In any case, the reader is referred to
the excellent review paper \cite{mclachlan02sm} and the monographs
\cite{hairer06gni,leimkuhler04shd}, where these and other issues,
mainly in connection with geometric numerical integration, are thoroughly
examined.

\section{Splitting and composition methods}

\subsection{Increasing the order of an integrator by composition}
  \label{s&c}

It is  well known that numerical integrators of arbitrarily high order can be obtained by composition of a basic integrator of low order. Consider for instance the leapfrog scheme (\ref{leapfrog}), which is a second-order integrator $\psis^{[2]}:\mathbb{R}^{2d} \to \mathbb{R}^{2d}$.
% that can be applied to integrate numericaly Hamiltonian systems of the form (\ref{eq.1.4}).
Then, a 4th order integrator $\psis^{[4]}:\mathbb{R}^{2d} \to \mathbb{R}^{2d}$ can be obtained as
\begin{equation}
\label{suzu1}
   \psis_{h}^{[4]} = \psis_{\alpha h}^{[2]} \circ \psis_{\beta h}^{[2]} \circ  \psis_{\alpha h}^{[2]}, \quad \mbox{ with } \quad
   \alpha = \frac{1}{2 - 2^{1/3}}, \qquad  \beta = 1 - 2 \alpha.
\end{equation}
More generally, if one recursively defines  $\psis^{[2k+2]}:\mathbb{R}^{2d} \to \mathbb{R}^{2d}$ for $k=1,2,\ldots$ as
\begin{equation}
  \label{eq:S_h-rec1}
   \psis_{h}^{[2k+2]} = \psis_{\alpha h}^{[2k]} \circ \psis_{\beta h}^{[2k]} \circ  \psis_{\alpha h}^{[2k]},
\end{equation}
with
\begin{equation}
  \label{eq:S_h-rec2}
     \alpha= \frac{1}{2 - 2^{1/(2k+1)}}, \qquad  \beta = 1 - 2 \alpha,
\end{equation}
then  the schemes $\psis_h^{[2k]}$ are of order $2k$ ($k\geq 1$) \cite{suzuki90fdo,yoshida90coh}.
We will prove later on this assertion.
At this point it is useful to introduce the
notion of adjoint of a given integrator $\psi_h$ \cite{sanz-serna94nhp}. By
definition, this is the method $\psi_h^{\ast}$ such that
$\psi_h^{\ast} = \psi_{-h}^{-1}$. A method that is equal to its own adjoint
is called self-adjoint or (time-)\textit{symmetric}. In this case,
$\psi_{-h} \circ \, \psi_h = \mathrm{id}$.
Since the leapfrog scheme (\ref{leapfrog}) can be rewritten as
\begin{equation}
  \label{eq:S_h}
\psis_h^{[2]} = \chi_{h/2} \circ \chi_{h/2}^*,
\end{equation}
where $\chi_h$ is the symplectic Euler method (\ref{Euler-sympl}), then $\psis_h^{[2]}$ is certainly time-symmetric. Actually, given any basic first order integrator $\chi_h:\mathbb{R}^D \to \mathbb{R}^D$ for the ODE system (\ref{eq.1.1}), the composition (\ref{eq:S_h}) is a time-symmetric method of order $2$, and the other way around:
any  self-adjoint second order integrator can be expressed as (\ref{eq:S_h}).
Furthermore, the integrators $\psis_h^{[2k]}$ ($k=1,2,\ldots$) recursively defined by (\ref{eq:S_h-rec1})--(\ref{eq:S_h-rec2}) are time-symmetric methods of order $2k$. In particular, if  $f(x)$ in the ODE (\ref{eq.1.1}) is split as
\begin{equation}
  \label{eq:split-f}
f(x) = \sum_{i=1}^m f^{[i]}(x)
\end{equation}
then, time-symmetric integrators $\psis_h^{[2k]}$ of order $2k$ can be constructed in this way by considering the basic first order integrator
\begin{equation}
\label{eq:chih}
   \chi_h = \varphi_h^{[m]} \circ \cdots \circ \varphi_h^{[2]} \circ \varphi_h^{[1]}
\end{equation}
and its adjoint
\begin{equation*}
   \chi^*_h = \chi_{-h}^{-1} = \varphi_h^{[1]} \circ \varphi_h^{[2]} \circ
\cdots \circ \varphi_h^{[m]}.
\end{equation*}
%
%Recall that, if the split ODE systems (\ref{eq.1.2}) share  with the original system (\ref{eq.1.1})  some qualitative features (such as being Hamiltonian systems, or divergence-free systems) which induce certain geometric properties of their flows (symplecticity, preservation of volume), then the high order integrators $\psis_h^{[2k]}$ also will keep the same geometric properties.
%However, t
This general procedure of constructing geometric integrators of arbitrarily high order, although simple, presents  some drawbacks. In
particular,
the resulting methods require a large number of evaluations and usually have large truncation errors.

As we will show, efficient schemes can be built by considering more general composition integrators. First observe that the $(2k)$th order integrators $\psis_h^{[2k]}$ can be rewritten in the form
\begin{equation}
\label{eq:compint}
\psi_h =  \chi_{\alpha_{2s} h}\circ
\chi^*_{\alpha_{2s-1}h}\circ
\cdots\circ
\chi_{\alpha_{2}h}\circ
\chi^*_{\alpha_{1}h}
\end{equation}
with $s=3^{k-1}$ and some fixed coefficients $(\alpha_1,\ldots,\alpha_{2s}) \in \r^{2s}$. Then the
 idea is to consider composition integrators of the form (\ref{eq:compint})
with appropriately chosen coefficients $(\alpha_1,\ldots,\alpha_{2s}) \in \r^{2s}$.

In the particular case  where the ODE (\ref{eq.1.1}) is split in two parts $f=f^{[a]}+f^{[b]}$ and $\chi_h =  \varphi^{[b]}_{h}\circ \varphi^{[a]}_{h}$,
one can trivially check that the composition integrator (\ref{eq:compint}) can be rewritten as
\begin{equation}
  \label{eq:splitting}
\psi_h = \varphi^{[b]}_{b_{s+1} h}\circ
 \varphi^{[a]}_{a_{s}h}\circ \varphi^{[b]}_{b_{s} h}\circ
 \cdots\circ
 \varphi^{[b]}_{b_{2}h}\circ
 \varphi^{[a]}_{a_{1}h} \circ \varphi^{[b]}_{b_{1}h},
\end{equation}
where $b_{1}=\alpha_{1}$ and for $j=1,\ldots,s$,
\begin{equation}
  \label{eq:abalpha}
a_j=\alpha_{2j-1}+\alpha_{2j}, \qquad\quad  b_{j+1} = \alpha_{2j} + \alpha_{2j+1}
\end{equation}
 (with $\alpha_{2s+1}=0$). Conversely, any integrator of the form (\ref{eq:splitting}) satisfying that $\sum_{i=1}^{s} a_i = \sum_{i=1}^{s+1} b_i$ can be expressed in the form (\ref{eq:compint}) with $\chi_h =  \varphi^{[b]}_{h}\circ \varphi^{[a]}_{h}$.  For later reference, we state the following result, due to McLachlan~\cite{mclachlan95otn}.

\begin{theorem}
\label{th:McLachlan}
 The integrator (\ref{eq:splitting}) is of order $r$ for ODEs of the form (\ref{eq.1.1}) with $f: \mathbb{R}^D \longrightarrow  \mathbb{R}^D$ arbitrarily split as $f=f^{[a]}+f^{[b]}$ if and only if the integrator (\ref{eq:compint}) (with coefficients $\alpha_j$ obtained from (\ref{eq:abalpha})) is of order $r$ for
arbitrary consistent integrators $\chi_h$.
\end{theorem}

\subsection{Integrators and series of differential operators}
\label{ss:diffop}

Before proceeding further with the analysis, let us relate generic numerical integrators with
formal series of differential equations. This relationship will allow one to formulate in a rather simple way
the conditions
to be satisfied by an integration scheme to achieve a given order of consistency.

First of all, let us recall that an
integrator $\psi_h:\mathbb{R}^{D} \to \mathbb{R}^{D}$ for the system (\ref{eq.1.1}) is said to be of order $r$ if for all $x \in \mathbb{R}^{D}$
\begin{equation}   \label{eq.2.1}
   \psi_h(x) = \varphi_h(x) + \cO(h^{r+1})
\end{equation}
as $h\rightarrow 0$, where $\varphi_h$
is the $h$-flow of the ODE (\ref{eq.1.1}).

It is well known that, for any smooth function $g:\mathbb{R}^D\rightarrow \mathbb{R}$, it formally holds that \cite{olver93aol}
\begin{eqnarray*}
g(\varphi_h(x))=g(x) + \sum_{n\geq 1} \frac{1}{n!} F^n[g](x) = \exp(h F)[g](x),
\end{eqnarray*}
where $F$ is the Lie derivative associated to the ODE system (\ref{eq.1.1}), i.e., the first order linear differential operator $F$ acting on functions in $\cinf$ as follows. For each $g \in \cinf$ and each
$x = (x_1, \ldots, x_D) \in \r^D$
\begin{eqnarray}
\label{eq:F}
F[g](x) = \sum_{j=1}^{D} f_j(x) \frac{\partial g}{\partial x_j}(x),
\end{eqnarray}
where $f(x)=(f_1(x),\ldots,f_D(x))^T$.
%
%As for the basic splitting integrator (\ref{eq:chih}), we have that
%\begin{eqnarray*}
%g(\chi_h(x))=\exp(h F^{[1]}) \cdots \exp(h F^{[m]})[g](x),
%\end{eqnarray*}
%where each $F^{[i]}$ ($1\leq i \leq m$) is the Lie derivative associated to the ODE $x'=f^{[i]}(x)$.
Motivated by this fact, we consider for a basic integrator $\chi_h:\mathbb{R}^D \to \mathbb{R}^D$,  the linear differential operators $X_n$ ($n\geq 1$) acting on smooth functions $g\in\cinf$ as follows:
\begin{equation}
\label{eq:X_n}
X_n[g](x) =  \frac{1}{n!}\frac{d^n}{dh^n} \left. g(\chi_h(x)) \right|_{h=0},
\end{equation}
so that formally $g(\chi_h(x)) = X(h)[g](x)$, where
% $X(h)$ is the series of differential operators
\begin{equation}
\label{eq:X(h)}
X(h) = I + \sum_{n\geq 1} h^n X_n,
\end{equation}
and $I$ denotes the identity operator.
Thus, the integrator $\chi_h$ is of order $r$ if
\begin{equation*}
X_n = \frac{1}{n!} F^n, \qquad 1 \leq n \leq r.
\end{equation*}
Alternatively, one can consider the series of vector fields
\[
Y(h) = \sum_{n\geq 1} h^n Y_n = \log(X(h)) =
\sum_{m\geq1} \frac{(-1)^{m+1}}{m} \left( h X_1 + h^2 X_2+\cdots \right)^m,
\]
that is,
\begin{equation*}
Y_n = \sum_{m\geq1}^n \frac{(-1)^{m+1}}{m} \sum_{j_1+\cdots+j_m=n}
X_{j_1} \cdots X_{j_m},
\end{equation*}
so that $X(h)= \exp(Y(h))$, and formally, $g(\chi_h(x)) =
\exp(Y(h))[g](x)$. Clearly,  the basic integrator is of order $r$
if
%It can be shown that the algebraic properties of the linear
%differential operators $\Psi_k$ imply that such
%$F_h$ is a series of vector fields. This means the well known
%fact that \cite{HLW} the integrator
%$\psi_h$ can be formally interpreted as the exact 1-flow of the
%modified vector field $\hf_h$. Then, condition (\ref{eq.2.2.1})
%is equivalent to
%
\begin{equation*}
 Y_1 = F, \qquad Y_n = 0 \quad  \mbox{for} \quad 2 \leq n \leq r.
\end{equation*}
For the adjoint integrator $\chi^*_h=\chi_{-h}^{-1}$, one obviously gets
 $g(\chi^*_h(x)) = \e^{-Y(-h)}[g](x)$. This shows that $\chi_h$ is time-symmetric if and only if $Y(h) = h Y_1 + h^3 Y_3 + \cdots$, and in particular, that time-symmetric methods are of even order.

It is possible now to check that the symmetric integrators $\psis_h^{[2k]}$ given by
(\ref{eq:S_h-rec1})--(\ref{eq:S_h-rec2})
are actually schemes of order $2k$ provided that $\psis_h^{[2]}$ is a symmetric second order integrator.  Consider the series of differential operators
\[
  F^{[2k]}(h)=h F + h^{2k+1} F_{2k+1}^{[2k]} + h^{2k+3} F_{2k+3}^{[2k]} + \cdots
\]
 such that $g(\psis_h^{[2k]}(x)) = \exp(F^{[2k]}(h))[g](x)$. Then one clearly has
  \[
    \exp(F^{[2k+2]}(h)) =     \exp(F^{[2k]}(\alpha h)) \exp(F^{[2k]}(\beta h)) \exp(F^{[2k]}(\alpha h))
  \]
which implies
\[
  F^{[2k+2]}(h) =   h \, (2\alpha+\beta) \, F +   h^{2k+1} \, (2\alpha^{2k+1}+\beta^{2k+1}) \, F_{2k+1}^{[2k]} + \mathcal{O}(h^{2k+3}),
\]
and thus $\psis_h^{2k+2}$ is of order $2k+2$ provided that  $\psis_h^{2k}$ is of
 order $2k$ and $\alpha$ and $\beta$ satisfy the equations
\[
  2\alpha+\beta=1, \qquad 2\alpha^{2k+1}+\beta^{2k+1}=0,
\]
whose unique real solution is given by (\ref{eq:S_h-rec2}).

In the general case, for the composition method (\ref{eq:compint}) we have
\begin{eqnarray*}
g(\psi_h(x))= \Psi(h)[g](x),
\end{eqnarray*}
where $\Psi(h)=I + h \Psi_1 + h^2 \Psi_2 + \cdots$ is a series of differential operators satisfying
\begin{equation}
\label{eq:Psi(h)}
\Psi(h)=  X(-\alpha_1 h)^{-1}X(\alpha_2 h) \cdots X(-\alpha_{2s-1} h)^{-1}X(\alpha_{2s} h),
\end{equation}
where the series $X(h)$ is given by (\ref{eq:X_n})--(\ref{eq:X(h)}), and
\begin{eqnarray}
  \label{eq:X(h)^-1}
X(h)^{-1} = I + \sum_{m\geq1} (-1)^{m+1} \left( h X_1 + h^2 X_2+\cdots \right)^m.
\end{eqnarray}
Thus, the order of a composition integrator of the form (\ref{eq:compint}) can be checked by comparing the series of differential operators $\Psi(h)$ with the series $\exp(h F)$ associated to the flow of the system (\ref{eq.1.1}). That is, the integrator (\ref{eq:compint}) is of order $r$ if
\begin{equation}  \label{eq.2.2.1}
\Psi_n = \frac{1}{n!} F^n, \qquad 1 \leq n \leq r.
\end{equation}
Instead of using (\ref{eq:Psi(h)}) to obtain the terms $\Psi_n$ of the series $\Psi(h)$, one can equivalently consider the formal equality
\begin{equation}   \label{eq.211b}
\Psi(h) = \e^{-Y(-h\alpha_1)} \, \e^{Y(h \alpha_2)} \cdots \, 
\e^{-Y(-h\alpha_{2s-1})} \, \e^{Y(h \alpha_{2s})},
\end{equation}
to obtain the series expansion of $\log(\Psi(h))=\sum_{n\geq 1} h^n F_n$, so that $r$th order compositions methods can also be characterized by the conditions
\begin{equation}  \label{eq.2.2.3}
 F_1 = F, \qquad F_n = 0 \quad  \mbox{for} \quad 2 \leq n \leq r.
\end{equation}

As for the splitting integrator (\ref{eq:splitting}) when the ODE (\ref{eq.1.1}) is split in two parts,
\begin{equation}
  \label{eq:f=fa+fb}
f(x)=f^{[a]}(x)+f^{[b]}(x),
\end{equation}
let $F^{[a]}$ and $F^{[b]}$ be the Lie derivatives corresponding to $f^{[a]}$ and $f^{[b]}$ respectively, that is,
\begin{eqnarray}
\label{eq:AB}
F^{[a]}[g](x) = \sum_{j=1}^{D} f_j^{[a]}(x) \frac{\partial g}{\partial x_j}(x),
\qquad
F^{[b]}[g](x) = \sum_{j=1}^{D} f_j^{[b]}(x) \frac{\partial g}{\partial x_j}(x)
\end{eqnarray}
for each $g \in \cinf$ and each $x \in \r^D$. Then, the series $\Psi(h)$ of differential operators associated to the integrator $\psi_h$ in (\ref{eq:splitting}) can be formally written as
\begin{equation}
  \label{eq:Psi(h)splitting}
\Psi(h) = \e^{b_{1} h F^{[b]}} \, \e^{a_{1}h F^{[a]}} \cdots \, 
\e^{b_{s}h F^{[b]}} \, \e^{a_{s}h F^{[a]}} \, \e^{b_{s+1}h F^{[b]}}.
\end{equation}

\section{Order conditions of splitting and composition methods}

There are several procedures to get the order conditions for the
coefficients of splitting and composition methods of a given
order. These are, generally speaking, large systems of polynomial
equations in the coefficients which are obtained from equations
(\ref{eq.2.2.3}). Perhaps the two most popular are (i) the
expansion of the series $\log(\Psi(h))=\sum_{n\geq 1} h^n F_n$ of
vector fields by applying recursively the
Baker--Campbell--Hausdorff (BCH) formula \cite{yoshida90coh}, and
(ii) a generalization of the theory of rooted trees used in the
theory of Runge--Kutta methods, which allows one to get an
equivalent set of simpler order conditions in a systematic
way~\cite{murua99ocf} (see also \cite{hairer06gni}). In this
section we first summarize briefly how to get these equations with
the BCH formula, and then we present a novel approach, related to
that in~\cite{murua99ocf}, but based on Lyndon words instead of
rooted trees.

\subsection{Order conditions  via BCH formula}  \label{sec.BCH}

As is well known, if $X$ and $Y$ are two non-commuting operators, the BCH formula establishes that formally,
$\e^X \e^Y = \e^Z$, where $Z$  belongs to the Lie algebra $\cL(X,Y)$ generated by $X$ and $Y$ with the commutator $[X,Y]=X Y - Y X$ as Lie bracket \cite{varadarajan84lgl}. Moreover,
\begin{equation}  \label{BCH1}
   Z = \log( \e^X \, \e^Y) = X + Y +
    \sum_{m=2}^{\infty} \, Z_m,
\end{equation}
with $Z_m(X,Y)$ a homogeneous Lie polynomial in $X$ and $Y$
of degree $m$, i.e., it is a $\mathbb{Q}$-linear combination
of commutators of the form $[V_1,[V_2, \ldots,[V_{m-1},V_m]
\ldots]]$ with $V_i \in \{X,Y\}$ for $1 \le i \le m$. The first terms read
\begin{eqnarray*}
  Z_2 & = & \frac{1}{2} [X,Y]   \\
  Z_3 & = & -\frac{1}{12} [[X,Y],X] + \frac{1}{12}
                [[X,Y],Y]    \\
  Z_4 & = & \frac{1}{24} [[[X,Y],Y],X],
\end{eqnarray*}
and explicit expressions up to $m=20$ have been recently computed in an arbitrary generalized Hall basis
of $\cL(X,Y)$ \cite{casas08aea}.

The procedure to get the order conditions for the composition method (\ref{eq:compint})
with this approach
can be summarized as follows. First, consider the
series of differential operators $\Psi(h)$ associated to the integrator (\ref{eq:compint}), expressed as a product of exponentials of vector fields, i.e., equation
(\ref{eq.211b}). Then, apply repeatedly the BCH formula to get the series expansion $\log(\Psi(h))=\sum_{n\geq 1} h^n F_n$.  In this way, one gets
\begin{eqnarray}  \label{bch2}
  \log( \Psi(h)) &=&  h w_{1} Y_1 + h^2 w_{2} Y_2
 + h^3 (w_{3} Y_3 +  w_{12} [Y_1, Y_2]) \\ \nonumber
&& + h^4 (w_{4} Y_4 + w_{1 3} [Y_1,Y_3] + w_{112} [Y_1,[Y_1,Y_2]]) +   \cO(h^5)
\end{eqnarray}
where the $w_{j_1 \cdots j_m}$ are polynomials of degree $n=j_1+\cdots+j_m$ in the parameters $\alpha_1,\ldots,\alpha_{2s}$.
The first such polynomials are
\begin{equation}   \label{bch3}
   w_{1} = \sum_{i=1}^{2s} \alpha_i, \qquad w_{2} = \sum_{i=1}^{2s} (-1)^{i} \alpha_i^2, \qquad
   w_{3} = \sum_{i=1}^{2s} \alpha_i^3.
\end{equation}
In general, the expressions of the polynomials $w_{j_1 \cdots j_m}$ in (\ref{bch2}) obtained
by repeated application of the BCH formula are rather cumbersome.
%, one typically gets cumbersome expressions that we do not reproduce here.

The order conditions for the composition integrator (\ref{eq:compint}) are then obtained by imposing 
equations (\ref{eq.2.2.3}) to guarantee that the scheme has order $r \ge 1$. Thus, the order conditions 
 are $w_{1} = 1$, and $w_{j_1 \cdots j_m} = 0$ whenever $2 \leq j_1+\cdots+j_m \leq r$.

One can proceed similarly to get the order conditions of the splitting scheme (\ref{eq:splitting}) in terms of the coefficients $a_i,b_i$: Consider
the series of differential operators $\Psi(h)$ associated to the integrator (\ref{eq:splitting}) expressed as (\ref{eq:Psi(h)splitting}); then, apply repeatedly the BCH formula to get the series expansion $\log(\Psi(h))=\sum_{n\geq 1} h^n F_n$, so that the order conditions will be obtained by imposing equations (\ref{eq.2.2.3}) to guarantee order $r \ge 1$. It can be seen that
the following holds for $\log(\Psi(h))$,
\begin{eqnarray}  \nonumber
  \log( \Psi(h)) &=&  h (v_{a} F^{[a]} + v_{b} F^{[b]}) + h^2 v_{ab} F^{[ab]} + h^3 (v_{abb} F^{[abb]} + v_{aba} F^{[aba]}) \\  \label{eq:bch-splitting}
&& + h^4 (v_{abbb} F^{[abbb]} + v_{abba} F^{[abba]} + v_{abaa} F^{[abaa]}) +   \cO(h^5),
\end{eqnarray}
where
\begin{eqnarray*}
F^{[ab]} &\!=\!&  [F^{[a]}, F^{[b]}], \quad
F^{[abb]} = [F^{[ab]},F^{[b]}], \quad
F^{[aba]} = [F^{[ab]},F^{[a]}], \\
F^{[abbb]} &\!=\!&  [F^{[abb]}, F^{[b]}], \quad
F^{[abba]} = [F^{[abb]},F^{[a]}], \quad
F^{[abaa]} = [F^{[aba]},F^{[a]}],
\end{eqnarray*}
and $v_{a},v_{b},v_{ab},v_{abb},v_{aba},v_{abbb},\ldots$ are polynomials in the parameters $a_i,b_i$ of the splitting scheme (\ref{eq:splitting}).
In particular, one gets
\begin{eqnarray}  \label{eq:bch-splitting2}
   v_{a} &\!=\!& \sum_{i=1}^{s} a_i, \qquad v_{b} = \sum_{i=1}^{s+1} b_i, \qquad
   v_{ab} = \frac12 v_{a} v_{b} - \sum_{1\leq i \leq j \leq s} b_i a_j, \\
 \nonumber
   2 v_{aba} &\!=\!& -\frac16 v_{a}^2 v_{b} +\sum_{1\leq i < j \leq k \leq s} a_i b_j a_k,  \qquad
   2 v_{abb} = \frac{1}{6} v_{a} v_{b}^2 - \sum_{1\leq i\leq j < k \leq s+1} b_i a_j b_k.
\end{eqnarray}

From (\ref{eq:bch-splitting}), we see that a characterization of the
order of the splitting scheme (\ref{eq:splitting}) can be obtained by
considering $v_a=v_b=1$ and $v_{ab}=v_{abb}=v_{aba}=\cdots =0$ up
to polynomials of that form of the required order. The set of order
conditions thus obtained will be independent in the general case if
the vector fields $F^{[a]},F^{[b]},F^{[ab]},F^{[abb]}, F^{[aba]}$
considered in (\ref{eq:bch-splitting}) correspond to a basis of the
free Lie algebra on the alphabet $\{a,b\}$. Notice that in
(\ref{eq:bch-splitting}) we have considered a Hall basis (the
classical basis of P. Hall) associated to the Hall words
$a,b,ab,abb,aba,abbb,abba,abaa,\cdots$
\cite{reutenauer93fla}.  The coefficients $v_{w}$ in
(\ref{eq:bch-splitting}) corresponding to each Hall word $w$ can be
systematically obtained using the results in~\cite{murua06tha} in
terms of rooted trees and iterated integrals.  An efficient algorithm
(based on the results in~\cite{murua06tha}) of the BCH formula and
related calculations that allows one to obtain
(\ref{eq:bch-splitting}) up to terms of arbitrarily high degree is
presented in~\cite{casas08aea}.

\subsection{A set of independent order conditions}
\label{ss:oc}

We next present a set of order conditions for composition integrators (\ref{eq:compint}) derived in \cite{chartier08aat}.

From (\ref{eq:X(h)})--(\ref{eq:X(h)^-1}), it follows that
\begin{equation}
\label{eq:Psi(h)2}
  \Psi(h) = I + \sum_{n\geq 1} h^n \sum_{j_1+\cdots+j_r=n} u_{j_1 \cdots j_r}(\alpha_1,\ldots, \alpha_{2s}) \, X_{j_1} \cdots X_{j_r},
\end{equation}
for some polynomial functions $u_{j_1 \cdots j_r}$ of the parameters $\alpha_1,\ldots, \alpha_{2s}$ of the method.
We next introduce some notation in order to explicitly write these polynomials.
For each positive integer $j$, we write $j^*=j-1$ if $j$ is even, and $j^*=j$ if $j$ is odd.
 Finally, for each pair $(i,j)$ of positive integers, we write $\alpha_j^{(i)} = (-1)^{j(i-1)} (\alpha_j)^i$. That is, $\alpha_j^{(i)}=(\alpha_j)^{i}$ if $j$ is even or $i$ is odd, and $\alpha_j^{(i)}=-(\alpha_j)^{i}$ if $j$ is odd and $i$ is even.
Now, it is not difficult to check that, for each multi-index $(i_1,\ldots,i_m)$ of length $m\geq 1$ and $(\alpha_1,\ldots,\alpha_{2s}) \in \r^{2s}$,
\begin{equation}
  \label{eq:uii}
  u_{i_1 \cdots i_m}(\alpha_1,\ldots, \alpha_{2s}) = \sum_{1\leq j_1\leq j_1^*\leq j_2 \leq \cdots \leq j_{m-1} \leq j_{m-1}^* \leq  j_m \leq 2s} \alpha_{j_1}^{(i_1)} \cdots \alpha_{j_m}^{(i_m)}.
\end{equation}
Obviously, each $u_{i_1 \cdots i_m}$ can be seen as a real-valued function defined on the set
\begin{equation}
  \label{eq:G}
  \{(\alpha_1,\ldots,\alpha_{2s}) \in \r^{2s}\ : \ s\geq 1\}.
\end{equation}
Observe that each $u_{i_1 \cdots i_m}(\alpha_1,\ldots,\alpha_{2s})$ is a polynomial of degree $n=i_1+\cdots +i_m$ in the variables $\alpha_1,\ldots,\alpha_{2s}$.

Now, the order conditions of the composition scheme (\ref{eq:compint}) can be obtained by comparing the series (\ref{eq:Psi(h)2}) with $\exp(h F)$, that is, (\ref{eq.2.2.1}). Since $X_1=F$, as the basic integrator $\chi_h$ is assumed to be of order 1, we have that the method is of order $r$ if for each multi-index $(i_1, \ldots, i_m)$ with $i_1+ \cdots +i_m=n \leq r$,
\begin{equation}
\label{eq:oc-aux}
  u_{i_1 \cdots i_m}(\alpha_1,\ldots, \alpha_{2s}) = \left\{
  \begin{array}{lcc}
    \frac{1}{n!} &\mbox{ if } & (i_1, \ldots, i_m)=(1,\ldots,1), \\
    0 & \mbox{ otherwise.}    &
  \end{array}
\right.
\end{equation}
However, such order conditions are not independent. For instance, it can be checked that
\begin{eqnarray*}
  u_{11}=\frac{1}{2}(u_1^2 + u_2), \qquad
u_{21} = -u_{12} + u_3 + u_1 u_2, \qquad
  u_{111}=\frac{1}{6}u_1^3 + \frac12 u_{12} + \frac13 u_3,
\end{eqnarray*}
which implies that the order conditions (\ref{eq:oc-aux}) for $u_{11}$, $u_{12}$, $u_{111}$ are fulfilled provided that the conditions for $u_1,u_2,u_3,u_{12}$ hold.

A set of independent order conditions can be obtained as follows.
Consider the lexicographical order $<$ (i.e., the order used when
ordering words in the dictionary) on the set of multi-indices. A
multi-index $(i_1,\ldots,i_m)$ is a Lyndon multi-index if
$(i_1,\ldots,i_k) < (i_{k+1}, \ldots, i_m)$ for each $1 \leq k <
m$. For each $n\geq 1$, we denote as $L_n$ the set of functions
$u_{i_1 \cdots i_m}$ such that $(i_1,\ldots,i_m)$ is a Lyndon
multi-index satisfying that $i_1+\cdots +i_m=n$.
The first sets $L_n$ are
\begin{eqnarray*}
  L_1&=&\{u_1\}, \quad
  L_2 = \{u_2\}, \quad
  L_3 = \{u_{1 2},u_{3}\}, \quad
  L_4 = \{ u_{112}, u_{13},u_{4}\}, \\
  L_5 &=& \{u_{1112}, u_{113}, u_{122},  u_{14}, u_{23},  u_{5} \}.
\end{eqnarray*}
In particular, we have
\begin{eqnarray*}
   u_{1}(\alpha_1,\ldots,\alpha_{2s}) &=& \sum_{j=1}^{2s} \alpha_{j}, \\
   u_{2}(\alpha_1,\ldots,\alpha_{2s}) &=& \sum_{j=1}^{2s} (-1)^j \alpha_{j}^2, \\
   u_{3}(\alpha_1,\ldots,\alpha_{2s}) &=& \sum_{j=1}^{2s}  \alpha_{j}^3, \\
   u_{1 2}(\alpha_1,\ldots,\alpha_{2s}) &=& \sum_{j_2=1}^{2s} (-1)^{j_2} \alpha_{j_2} \sum_{j_1=1}^{j_2^*}  \alpha_{j_1}, \\
   u_{4}(\alpha_1,\ldots,\alpha_{2s}) &=& \sum_{j=1}^{2s} (-1)^j \alpha_{j}^4, \\
   u_{1 3}(\alpha_1,\ldots,\alpha_{2s}) &=& \sum_{j_2=1}^{2s} \alpha_{j_2}^3 \sum_{j_1=1}^{j_2^*}  \alpha_{j_1}, \\
   u_{1 1 2}(\alpha_1,\ldots,\alpha_{2s}) &=&
\sum_{j_3=1}^{2s}  (-1)^{j_3} \alpha_{j_3}^2
\sum_{j_2=1}^{j_3^*} \alpha_{j_2} \sum_{j_1}^{j_2^*}  \alpha_{j_1},
\end{eqnarray*}
and so on.

We can finally state the following result~\cite{chartier08aat}: Given $(\alpha_1,\ldots,\alpha_{2s})$, the integrator (\ref{eq:compint}) is of order $r$ for arbitrary ODE systems (\ref{eq.1.1}) and arbitrary consistent integrators $\chi_h$ if and only if
$\alpha_1+\cdots+\alpha_{2s}=1$ (i.e. $u_1(\alpha_1,\ldots,\alpha_{2s})=1$) and
\begin{equation}
  \label{eq:occomp}
    \forall \, u \in \bigcup_{n\geq 2}^{r} L_n, \qquad \quad u(\alpha_1,\ldots,\alpha_{2s})=0.
\end{equation}
Furthermore, such order conditions are independent to each other if arbitrary sequences $(\alpha_1,\ldots,\alpha_{2s})$ of coefficients of the method are considered.

\subsection{Order conditions of compositions methods with symmetry}

The order conditions are simplified for ($2s$)-tuplas $(\alpha_1,\ldots,\alpha_{2s})$ such that
\begin{equation}
  \label{eq:sym}
  \alpha_{2s-j+1}=\alpha_{j}, \quad \mbox{ for all }  j.
\end{equation}
It is easy to check that the simplifying assumption (\ref{eq:sym}) implies that the composition integrator (\ref{eq:compint}) is time-symmetric (i.e., $\psi^*_h=\psi_h$). In that case,  only the conditions for $u \in L_{n}$ with odd $n$ remain independent.

 The order conditions can be alternatively simplified by requiring that
\begin{equation}
  \label{eq:sym2}
  \alpha_{2j}=\alpha_{2j-1}, \quad \forall \, j,
\end{equation}
in which case, only the conditions for Lyndon multi-indices $(i_1,\ldots,i_m)$ with odd $i_1,\ldots,i_m$ are required. The simplifying assumption (\ref{eq:sym2}) means that the composition integrator (\ref{eq:compint}) can be rewritten as
\begin{equation}
\label{eq:compSS}
  \psi_h = \psis^{[2]}_{h \beta_s} \circ \cdots \circ \psis^{[2]}_{h \beta_1},
\end{equation}
where $\beta_j = 2 \alpha_{2j}$ and $\psis_h^{[2]}$ is the self-adjoint second order integrator
 $\psis_h^{[2]} = \chi_{h/2} \circ \chi^*_{h/2}$.

The order conditions are thus considerably reduced if one considers
composition methods satisfying both assumptions (\ref{eq:sym})--(\ref{eq:sym2}), that is, methods of the form (\ref{eq:compSS}) satisfying that
\begin{equation}
  \label{eq:sym3}
\beta_{j}=\beta_{s-j+1}, \quad \forall \, j.
\end{equation}
Schemes of this form can be dubbed as symmetric
 compositions of symmetric schemes.
For instance, for order $r\geq 6$ one has the conditions
\begin{eqnarray*}
 \sum_{j=1}^{2s} \alpha_{j}=1, \qquad
\sum_{j=1}^{2s}  \alpha_{j}^3 = 0, \qquad
\sum_{j=1}^{2s}  \alpha_{j}^5 = 0, \qquad
\sum_{j_3=1}^{2s}  \alpha_{j_3}^3
\sum_{j_2=1}^{j_3^*} \alpha_{j_2} \sum_{j_1}^{j_2^*}  \alpha_{j_1}=0
\end{eqnarray*}
in terms of the $\alpha_i$ coefficients (the actual expressions in terms of $\beta_i$ are slightly more involved).
In Table~\ref{tabla1} we display for each $k\geq 1$ the number $n_k$ of Lyndon multi-indices $(i_1,\ldots,i_m)$ with $i_1+\cdots+i_m=k$, and the number $m_k$ of Lyndon multi-indices $(i_1,\ldots,i_m)$ with $i_1+\cdots+i_m=k$ and odd indices $i_1,\ldots,i_m$. Thus, the number of independent conditions to guarantee that the general composition integrator (\ref{eq:compint}) is at least of order $r$ is $n_1+\cdots+n_r$, while in the case of the composition (\ref{eq:compSS}) based on a symmetric second order integrator $\psis^{[2]}_h$ (or equivalently, a composition integrator (\ref{eq:compint}) with the additional symmetry condition (\ref{eq:sym2})), the number of independent order conditions is $m_1+\cdots+m_r$. If time-symmetry is imposed in the method (\ref{eq:compint}) (resp. (\ref{eq:compSS})) by the additional assumption (\ref{eq:sym}) (resp. (\ref{eq:sym3})), then there are $n_1+n_3+\cdots+n_{2l-1}$ (resp. $m_1+m_3+\cdots+m_{2l-1}$) independent conditions that guarantee order at least $r=2l$.
\begin{table}
\begin{center}
  \begin{tabular}{|c|ccccccccccc|} \hline
  $k$ & 1 & 2 & 3 & 4 & 5 & 6 & 7  &  8 & 9 & 10 & 11 \\ \hline
$n_k$ & 1 & 1 & 2 & 3 & 6 & 9 & 18 & 30 & 56& 99 & 186 \\
$m_k$ & 1 & 0 & 1 & 1 & 2 & 2 &  4 & 5  & 8 & 11 & 17 \\ \hline
\end{tabular}
\end{center}
  \caption{The numbers $n_k$ and $m_k$ of independent order conditions for general composition methods (\ref{eq:compint}) and for compositions (\ref{eq:compSS}) of a basic time-symmetric method, respectively.}
  \label{tabla1}
\end{table}

\subsection{Relation among different sets of order conditions of composition methods}

 In~\cite{murua99ocf}, a set of independent necessary and sufficient order conditions is given using labelled rooted trees (see also~\cite{hairer06gni}). A family of sets $\{\mathcal{T}_n\}_{n=1,2,\ldots}$ of functions defined  on the set (\ref{eq:G}) is identified such that the integrator (\ref{eq:compint}) is of order $r$ if and only if
$\alpha_1+\cdots+\alpha_{2s}=1$ and
\begin{equation}
\label{eq:occomp2}
    \forall \, u \in \bigcup_{n\geq 2}^{r} \mathcal{T}_n, \qquad\quad u(\alpha_1,\ldots,\alpha_{2s})=0.
\end{equation}
Each $u(\alpha_1,\ldots,\alpha_{2s})$ with $u \in \mathcal{T}_n$ is (as in the case where $u \in \mathcal{L}_n$), a polynomial of homogeneous degree $n$.
In particular,
\begin{eqnarray*}
  \mathcal{T}_1&=&\{u_1\}, \quad
\mathcal{T}_2 = \{u_2\}, \quad
\mathcal{T}_3 = \{u_{2 1},u_{3}\}, \quad
\mathcal{T}_4 = \{ v_{211}, u_{31},u_{4}\},
%\\ \mathcal{T}_5 &=& \{v_{2111}, v_{311}, v_{221},  u_{41}, u_{32},  u_{5} \},
\end{eqnarray*}
where the functions of the form $u_{i_1 \cdots i_m}$ are defined in (\ref{eq:uii}), and
\begin{eqnarray*}
v_{211}(\alpha_1,\ldots,\alpha_{2s}) &=&
\sum_{j_2=1}^{2s} (-1)^{j_2} \alpha_{j_2}^2 \left( \sum_{j_1=1}^{j_2^*} \alpha_{j_1}\right)^2.
%\\ v_{2111}(\alpha_1,\ldots,\alpha_{2s}) &=&
%\sum_{j_2=1}^{2s} (-1)^{j_2} \alpha_{j_2}^2 \left( \sum_{j_1=1}^{j_2^*} \alpha_{j_1}\right)^3, \\
%v_{311}(\alpha_1,\ldots,\alpha_{2s}) &=&
%\sum_{j_2=1}^{2s} \alpha_{j_2}^3 \left( \sum_{j_1=1}^{j_2^*} \alpha_{j_1}\right)^2, \\
%v_{221}(\alpha_1,\ldots,\alpha_{2s}) &=&
%\sum_{j_3=1}^{2s} (-1)^{j_3} \alpha_{j_3}^2 \sum_{j_1,j_2=1}^{j_3^*} (-1)^{j_2} \alpha_{j_1} \alpha_{j_2}^2.
\end{eqnarray*}
As shown in~\cite{chartier08aat}, the order conditions (\ref{eq:occomp2}) are equivalent to the conditions (\ref{eq:occomp}), as both $\cup_{n\geq 1} L_n$ and $\cup_{n\geq 1} \mathcal{T}_n$ generate the same graded algebra $\mathcal{H} = \bigoplus_{n\geq 1} \mathcal{H}_n$ of functions on the set (\ref{eq:G}) (for each $u \in \mathcal{H}_n$, $u(\alpha_1,\ldots,\alpha_{2s})$ is a polynomial of homogeneous degree $n$, actually, a linear combination of polynomials $u_{i_1\cdots i_m}$ of homegeneous degre $n=i_1+\cdots +i_m$). For instance, it can be seen that
\begin{equation*}
  v_{211} = 2 u_{211} - u_{22} = 2 (u_{112} - u_{13} +  u_1 u_{12} +u_3 u_1) + u_1^2 u_2 + \frac12 (u_4-u_2^2).
\end{equation*}
Finding an independent set of order conditions for composition integrators is equivalent to finding a set of functions of homogeneous degree that generate the algebra $\mathcal{H}$ (see~\cite{chartier08aat}) for more details.

Of course, the functions $w_{i_1 \cdots i_m}$ in (\ref{bch2}) obtained when deriving the order conditions of composition integrators by repeated use of the BCH formula also belong to the same algebra of functions. For instance, $w_{n}=u_{n}$, and $w_{12} = u_{12}-u_3 - u_1 u_2$.

Recall that Theorem~\ref{th:McLachlan} characterizes the order conditions of splitting integrators of the form (\ref{eq:splitting}), where the ODE (\ref{eq.1.1}) is split in two parts (\ref{eq:f=fa+fb}), in terms of the order conditions of composition integrators (\ref{eq:compint}). Actually, the polynomials  $v_{a},v_{b},v_{ba},v_{baa},v_{bab},v_{baaa},\ldots$ (on the parameters $a_i,b_i$) in (\ref{eq:bch-splitting}) can be rewritten as linear combinations of the polynomials (on the parameters $\alpha_i$) in (\ref{eq:uii}) provided that  (\ref{eq:abalpha}) and $v_{a}=v_{b}=u_1$ hold. In particular, it can be seen that
\begin{eqnarray*}
v_{ab} &=&  \frac{u_2}{2}, \\
v_{abb} &=& \frac{1}{12} \left(-u_3-3 u_{12}+3
   u_{21}\right), \\
v_{aba} &=& \frac{1}{12} \left(u_3-3 u_{12}+3
   u_{21}\right), \\
v_{abbb} &=& \frac{1}{12}
   \left(u_{22}-u_{31}+u_{112}-2
   u_{121}+u_{211}\right), \\
v_{abba} &=& \frac{1}{24} \left(-u_4-2 u_{13}+4
   u_{22}-2 u_{31}+4 u_{112}-8
   u_{121}+4 u_{211}\right), \\
v_{abaa} &=& \frac{1}{12}
   \left(-u_{13}+u_{22}+u_{112}-2
   u_{121}+u_{211}\right).
%\\ \frac{1}{720} \left(-u_5-5 u_{14}-10
%   u_{23}+20 u_{32}-5 u_{41}-20
%   u_{113}+30 u_{122}-20 u_{131}-30
%   u_{221}+40 u_{311}-60 u_{1121}+60
%   u_{1211}\right) \\
% \frac{1}{480} \left(-u_5-25
%   u_{14}-40 u_{23}+50 u_{32}+25
%   u_{41}-110 u_{113}+30 u_{122}-10
%   u_{131}-30 u_{221}+130 u_{311}-90
%   u_{1112}-150 u_{1121}+150
%   u_{1211}+90 u_{2111}\right) \\
% \frac{1}{480} \left(u_5-25 u_{14}-50
%   u_{23}+40 u_{32}+25 u_{41}-130
%   u_{113}+30 u_{122}+10 u_{131}-30
%   u_{221}+110 u_{311}-90
%   u_{1112}-150 u_{1121}+150
%   u_{1211}+90 u_{2111}\right) \\
% \frac{1}{720} \left(u_5+5 u_{14}-20
%   u_{23}+10 u_{32}+5 u_{41}-40
%   u_{113}+30 u_{122}+20 u_{131}-30
%   u_{221}+20 u_{311}-60 u_{1121}+60
%   u_{1211}\right) \\
% \frac{1}{240} \left(-2 u_5-10
%   u_{14}-5 u_{23}+10 u_{32}+5
%   u_{41}-25 u_{113}+5 u_{131}+15
%   u_{212}-30 u_{221}+35 u_{311}-15
%   u_{1112}-45 u_{1121}+45
%   u_{1211}+15 u_{2111}\right) \\
% \frac{1}{480} \left(-3 u_5-15
%   u_{14}-30 u_{23}+30 u_{32}+5
%   u_{41}-60 u_{113}-30 u_{122}+20
%   u_{131}+30 u_{212}-30 u_{221}+60
%   u_{311}-60 u_{1112}-60 u_{1121}+60
%   u_{1211}+60 u_{2111}\right)
\end{eqnarray*}

\subsection{Negative time steps}
\label{ntsteps}

It has been
noticed that some of the coefficients in splitting schemes (\ref{eq:splitting}) are
negative when the order $r \ge 3$. In other words, the methods always involve
stepping backwards in time. This constitutes a problem when the
differential equation is defined in a semigroup, as arises sometimes
in applications, since
then the method can only be conditionally stable \cite{mclachlan02sm}. Also
schemes with negative coefficients may not be well-posed when
applied to PDEs involving unbounded operators.

The existence of backward fractional time steps in this class of methods is
unavoidable, as shown in \cite{goldman96noo,sheng89slp,suzuki91gto}.  In fact, it can be established in an elementary way by virtue of the relationship
between the order conditions of schemes (\ref{eq:splitting}) and (\ref{eq:compint}) stated in Theorem~\ref{th:McLachlan} \cite{blanes05otn}:
Any splitting method of the form (\ref{eq:splitting}) that has order $r\geq 3$ neccesarily must fullfil the condition
\begin{equation}   \label{nts1}
u_{3}(\alpha_1,\ldots,\alpha_{2s}) =  \sum_{i=1}^{2s} \alpha_i^3 = \sum_{i=1}^s (\alpha_{2i-1}^3 + \alpha_{2i}^3) = 0,
\end{equation}
with coefficients $\alpha_j$ obtained from the relations (\ref{eq:abalpha}).
Since, for all $x, y \in \mathbb{R}$, it is true that $x^3 + y^3 < 0$ implies $x + y < 0$, then there must
exist some $i \in \{1, \ldots, s\}$ in the sum of (\ref{nts1}) such that
\[
  \alpha_{2i-1}^3 + \alpha_{2i}^3 < 0  \qquad \mbox{ and thus } \qquad \alpha_{2i-1} + \alpha_{2i} =
  a_i < 0.
\]
Obviously, one can also write (by taking $\alpha_0 = 0$)
\[
   u_{3}(\alpha_1,\ldots,\alpha_{2s}) = \sum_{i=0}^{2s+1}   \alpha_i^3 = \sum_{i=1}^{s+1} (\alpha_{2i-1}^3 + \alpha_{2i-2}^3) = 0
\]
just by grouping terms in a different way, and thus, by repeating the argument,
 there must exist some $j \in \{1, \ldots, s+1\}$ such that
\[
   \alpha_{2j-1} + \alpha_{2j-2} = b_j < 0.
\]
This proof shows clearly the origin of the existence of backward time steps: the equation $u_3=0$ can be
satisfied only if at least one $a_i$ and one $b_i$ are negative. According to this conclusion,
any splitting method of the form (\ref{eq:splitting}) verifying the order condition $u_{3}=0$ has necessarily
some negative coefficient $a_i$ and also some negative $b_i$.

\subsection{Near-integrable systems}
\label{Q&I}

In Hamiltonian dynamics one often encounters systems whose Hamiltonian function $H$ is a small perturbation of an exactly integrable Hamiltonian $H_0$, that is $H = H_0 + \varepsilon H_1$ with $\varepsilon \ll 1$. The perturbed Kepler problem (\ref{eq.14}) belongs to this category of near-integrable Hamiltonian systems. The gravitational $N$-body problem (\ref{eq.14N}),  when using
Jacobi coordinates, also falls within this class of problems. In that case, $H_0$ represents the Keplerian motion and $\varepsilon H_1$ the
mutual perturbations of the bodies on one another \cite{wisdom91smf}.

More generally, let us consider an ODE system
\begin{equation}     \label{qi.1}
      x^\prime = f^{[a]}(x) + \varepsilon f^{[b]}(x),
\end{equation}
containing a small parameter $|\varepsilon|\ll 1$. If the exact $h$-flows $\varphi^{[a]}_h$ and $\varphi^{[b]}_h$ of $x^{\prime}=f^{[a]}(x)$ and $x^{\prime}=\varepsilon\, f^{[b]}(x)$ respectively can be efficiently computed, then a scheme $\psi_h$ of the form
(\ref{eq:splitting}) can perform particularly well provided that the coefficients $a_i,b_i$ are appropriately chosen. To see this, consider the Lie derivatives (\ref{eq:AB}) of $f^{[a]}$ and $f^{[b]}$ respectively, so that the
corresponding series $\Psi(h)$  (\ref{eq:Psi(h)splitting}) of differential operators associated to the scheme (\ref{eq:splitting}) becomes
\begin{eqnarray*}
  \Psi(h) = \e^{b_{1} h \varepsilon F^{[b]}} \, \e^{a_{1}h F^{[a]}} \cdots \,
\e^{b_{s}h \varepsilon  F^{[b]}} \, \e^{a_{s}h F^{[a]}} \, \e^{b_{s+1}h \varepsilon  F^{[b]}}.
\end{eqnarray*}
Successive application of the BCH formula then leads to (\ref{eq:bch-splitting}) with $F^{[b]}$ replaced by $\varepsilon F^{[b]}$, that is
\begin{eqnarray*}
    \log( \Psi(h)) &=&
 h v_{a} F^{[a]} + \varepsilon ( h v_{b} F^{[b]} + h^2 v_{ab} F^{[ab]} + h^3 v_{aba} F^{[aba]} +  h^4 v_{abaa} F^{[abaa]}) \\
&& +  \varepsilon^2 (h^3  v_{abb} F^{[abb]} + h^4 v_{abba} F^{[abba]}) + \varepsilon^3
h^4 v_{abbb} F^{[abbb]} +   \cO(\varepsilon h^5).
\end{eqnarray*}
In practical applications one usually has $\varepsilon \ll h$ (or at least $\varepsilon \approx h$), so that
one is mainly interested in eliminating error terms with small powers of $\varepsilon$. For instance, if the coefficients $a_i,b_i$ of the splitting methods are chosen in such a way that
\begin{equation*}
  v_{a}=1=v_{b}, \quad v_{ab}=v_{aba}=v_{abaa}=v_{abb}=0,
\end{equation*}
then
\begin{equation*}
    \log( \Psi(h)) - h F = \cO(\varepsilon h^5 + \varepsilon^2 h^4),
\end{equation*}
where $F=F^{[a]} + \varepsilon F^{[b]}$.
More generally, one is interested in designing methods such that \cite{blanes00psm}
\begin{equation} \label{c.1}
   \log( \Psi(h)) - h F =  \cO(\varepsilon h^{s_{1}+1}+
  \varepsilon^{2} h^{s_{2}+1}+\varepsilon^{3}h^{s_{3}+1}+\cdots+\varepsilon^{m}
h^{s_{m}+1}).
\end{equation}
%A method which satisfies this condition
%is said to be of order $(s_{1},s_{2},s_{3},\ldots,s_{r}=r)$.
%We are interested in the case where
%$s_{i}\geq s_{i+1}$ and the list terminates with
%$\varepsilon^{r}h^{r+1}$,  $r$ being the standard order of consistency of the
%method.
Observe that $s_{1}$ is the order of consistency the
method would have in the limit $\varepsilon\rightarrow0$.
It is relatively easy to eliminate errors of order $\varepsilon h^k$ because there is only one such term for each order $k$,
namely $h^k \varepsilon \, v_{aba\cdots a} F^{[aba\cdots a]}$ (with $F^{[aba\cdots a]} = [[\cdots [[F^{[a]}, F^{[b]}],F^{[a]}] \ldots],F^{[a]}]$).

%For illustration, the $(4,2)$ methods are
%\begin{eqnarray*}
%  \psi_h & = & \varphi_{\alpha h}^{[1]} \circ \varphi_{\frac{1}{2} h}^{[2]} \circ \varphi_{(1-2\alpha) h}^{[1]} \circ
%      \varphi_{\frac{1}{2} h}^{[2]} \circ \varphi_{\alpha h}^{[1]}, \qquad  \alpha = (3-\sqrt{3})/6 \\
%   \psi_h & = &  \varphi_{\frac{1}{6} h}^{[2]} \circ \varphi_{\frac{1}{2} h}^{[1]} \circ \varphi_{\frac{2}{3} h}^{[2]} \circ
%          \varphi_{\frac{1}{2} h}^{[1]} \circ \varphi_{\frac{1}{6} h}^{[2]}
%\end{eqnarray*}

If one is interested in designing methods that approximate the exact solution up to
higher powers of $\varepsilon$, more terms have to be considered. In particular, there are
$\left\lfloor\frac12(k-1)\right\rfloor$ terms of order $\mathcal{O}(\varepsilon^2 h^k)$ and
$\left\lfloor\frac16(k-1)(k-2)\right\rfloor$ terms of order $\mathcal{O}(\varepsilon^3 h^k)$ \cite{mclachlan95cmi}.

\subsection{Runge--Kutta--Nystr\"om methods}
  \label{RKN}

Suppose now that one is interested in integrating numerically second-order ODE systems  of the form
\begin{equation}
    y^{\prime\prime}= g(y),  \label{rkn.1}
\end{equation}
where $y \in \mathbb{R}^{d}$ and $g: \mathbb{R}^{d}\longrightarrow \mathbb{R}^{d}$. In this case it is
still possible to use schemes (\ref{eq:splitting})  applied to the equivalent first-order ODE system.
More specifically, introducing the new variables $x=(y,v)$, with $v^\prime = y$,
and the maps $f^{[a]}:\mathbb{R}^{2d} \to \mathbb{R}^{2d}$ and $f^{[b]}:\mathbb{R}^{2d} \to \mathbb{R}^{2d}$ defined as
\begin{equation}
  \label{eq:f=fa+fbRKN}
 f^{[a]}(y,v) = (v,0), \qquad f^{[b]}(y,v) = (0,g(y)),
\end{equation}
equation (\ref{rkn.1}) can be rewritten as $x^\prime =
f^{[a]}(x) + f^{[b]}(x)$. Then, the splitting scheme (\ref{eq:splitting}) can be efficiently implemented, as the exact $h$-flows $\varphi^{[a]}_h$ and $\varphi^{[b]}_h$
of $x^\prime = f^{[a]}(x)$ and $x^\prime = f^{[b]}(x)$ are simply given by
\begin{equation}   \label{rkn.2}
   \begin{array}{ccl}
     \varphi_h^{[a]}(y,v) & = & (y + h v, v), \\
     \varphi_h^{[b]}(y,v) & = & (y, v + h g(y)).
   \end{array}
\end{equation}

It is not difficult to check that the splitting schemes of the
form (\ref{eq:splitting}) are particular instances of
Runge--Kutta--Nystr\"om (RKN) methods (see for
instance~\cite{hairer06gni}).

One of the most important applications of this class of schemes is
the study of Hamiltonian systems of the form $H(q,p) = T(p) +
V(q)$, where the kinetic energy $T(p)$ is quadratic in the momenta
$p$, i.e., $T(p)=\frac12 p^T M p$ for a symmetric square constant
matrix $M$, and $V(q)$ is the potential. In that case, the
corresponding Hamiltonian system can be written in the form
(\ref{rkn.1}) with $y=q$, $y' = v= M p$, and $g(y)=-\nabla V(y)$.

Although a splitting integrator (\ref{eq:splitting}) designed for
arbitrary ODE systems $x^\prime = f(x)$ split into two parts
(\ref{eq:f=fa+fb}) will perform well when applied to a second
order ODE system of the form (\ref{rkn.1}) with the splitting
(\ref{eq:f=fa+fbRKN}), much more efficient methods can be designed
in that case \cite{mclachlan92tao,calvo93tdo}. The main point here
is that in the case (\ref{eq:f=fa+fbRKN}) (we will refer to that
as the RKN case), $[[[F^{[b]},F^{[a]}],F^{[b]}],F^{[b]}] =  0$
identically. This is equivalent to $F^{[babb]}=0$ in
(\ref{eq:bch-splitting}), which introduces some linear
dependencies among higher order terms in the expansion of
$\log(\Psi(h))$ (see~\cite{mclachlan03tae} for a detailed study).
This means that the characterization given in
Theorem~\ref{th:McLachlan} for a splitting integrator
(\ref{eq:splitting}) to be of order $r$ (for $r\geq 4$) is no
longer applicable if one restricts to the case
(\ref{eq:f=fa+fbRKN}). In Table~\ref{tb:number-oc-RKN}, the number
of necessary and sufficient independent order conditions for a
splitting method (\ref{eq:splitting}) to be of order $r$ in the
RKN case is compared to the general case: For arbitrary systems
split into two parts (\ref{eq:f=fa+fb}), there are
$2+n_2+\cdots+n_r$ independent conditions (including the two
consistency conditions $v_a=v_b=1$), while in the RKN case, the
number of independent order conditions is $2+l_2+\cdots+l_r$.
Unfortunately, up to order three the order conditions are the same
in both cases and then, the results for negative time steps still
apply.
\begin{table}
\begin{center}
  \begin{tabular}{|c|cccccccccc|} \hline
  $k$ & 2 & 3 & 4 & 5 & 6 & 7  &  8 & 9 & 10 & 11 \\ \hline
$n_k$ & 1 & 2 & 3 & 6 & 9 & 18 & 30 & 56& 99 & 186 \\
$l_k$ & 1 & 2 & 2 & 4 & 5 & 10 & 14 & 25 & 39 & 69  \\ \hline
\end{tabular}
\end{center}
  \caption{The numbers $l_k$ of independent order conditions (at order $k$)
for splitting methods in the RKN case compared to the numbers
$n_k$ of conditions in the general case.} \label{tb:number-oc-RKN}
\end{table}

Since the reduction in the number of order conditions is due to
the fact that $f^{[a]}(y,v)$ is linear in $v$, it is immediate to
see that the methods obtained in this way also apply to the more
general problem 
$f^{[a]}(y,v)=(w_1(x)v+w_2(y),w_3(x)v+w_4(y))$, which includes the
system
\begin{equation}
    y^{\prime\prime}= M  y^{\prime} + g_1(y) + g_2(y).  \label{rkn.1b}
\end{equation}
Here splitting RKN methods are useful if the reduced problem
$y^{\prime\prime}= M  y^{\prime} + g_1(y)$
(i.e., $f^{[a]}(y,v)=(v,Mv+g_1(y))$) is easily solvable. For
Hamiltonian systems, this generalization corresponds to 
$H(q,p) = T(q,p) + V(q)$, where $T(q,p)=\frac12 p^T
M(q) p + f^T(q)p +W(q)$. Obviously, the exact solution for
$T(q,p)$ is only known for some particular cases, e.g. if $T$
corresponds to the Kepler problem in (\ref{eq.14}) or
(\ref{eq.14N}) or to the harmonic oscillator in
(\ref{HenonHeiles}) or (\ref{wave1}).

It is interesting to note that in quantum mechanics the kinetic
and potential energy verify analogue commutator rules to classical
mechanics and then RKN methods can also be used. If applied, for
instance to problems (\ref{Schrodinger}) and
(\ref{Gross-Pitaevskii}), one should keep in mind that in
the resulting composition method $\varphi_h^{[a]}$ must correspond
to the kinetic part.

We should also remark that, if the Hamiltonian function
$H(p,q)=\frac12 p^T M p + V(q)$ is such that in addition
$V(q)=\frac12 q^T N q$ (i.e., it corresponds to the generalized
harmonic oscillator (\ref{harmonicGeneral})), then the number of
order conditions reduces drastically: it is not difficult to see
that there is only one independent condition to increase the order
from $r=2k-1$ to $r=2k$, and only two to increase the order from
$r=2k$ to $r=2k+1$ (see~\cite{blanes07otl} for more details). We
will return later to this system and take profit of its special
features to design specially adapted splitting methods.

\section{Additional techniques to reduce the number of order conditions}

\subsection{Methods with modified potentials}
  \label{modified}

The splitting method (\ref{eq:splitting}) can be generalized by composing the exact flows of other vector fields in addition to $F^{[a]}$ and $F^{[b]}$, provided that they lie on the Lie algebra generated by $F^{[a]}$ and $F^{[b]}$. For instance, one could consider compositions that, in addition to $\varphi^{[a]}_{h}$ and $\varphi^{[b]}_{h}$, use the $h$-flow  $\varphi^{[abb]}_{h}$ of the vector field $F^{[abb]}=[[F^{[a]},F^{[b]}],F^{[b]}]$. To illustrate this fact, consider the composition
%\begin{equation}
%  \label{eq:modpotex}
%  \psi_h =
% \varphi_{b_1 h}^{[b]} \circ \varphi_{a_2 h}^{[a]} \circ  \varphi_{b_2 h}^{[b]} \circ \varphi_{c_2 h}^{[abb]} \circ \varphi_{b_2 h}^{[b]} \circ \varphi_{a_2 h}^{[a]} \circ \varphi_{b_1 h}^{[b]}
%\end{equation}
 \begin{equation}
  \label{eq:modpotex}
  \psi_h =
 \varphi_{h/6}^{[b]} \circ  \varphi_{h/2}^{[a]} \circ  \varphi_{h/3}^{[b]} \circ \varphi_{-h/72}^{[abb]} \circ \varphi_{h/3}^{[b]} \circ \varphi_{h/2}^{[a]} \circ \varphi_{h/6}^{[b]}.
\end{equation}
The scheme (\ref{eq:modpotex}), constructed in~\cite{chin97sif,koseleff94fcf}, is of order four. Indeed,
by repeated application of the BCH formula to 
\begin{eqnarray*}
\Psi(h) =  \e^{\frac{h}{6} F^{[b]}}  \e^{\frac{h}{2} F^{[a]}}  \e^{\frac{h}{3} F^{[b]}}  \e^{-\frac{h}{72} F^{[abb]}}  \e^{\frac{h}{3} F^{[b]}}  \e^{\frac{h}{2}F^{[a]}}  \e^{\frac{h}{6} F^{[b]}}
\end{eqnarray*}
one can check that $\Psi(h) = \e^{h(F^{[a]} + F^{[b]}) + \cO(h^5)}$.

Recall that, although the $h$-flows $\varphi^{[a]}_{h}$ and $\varphi^{[b]}_{h}$ of the vector fields $F^{[a]}$ and $F^{[b]}$ are by assumption computed easily, this is not necessarily the case for the $h$-flow $\varphi^{[abb]}_{h}$ of the vector field $F^{[abb]}$. However, in the RKN case (\ref{rkn.1}) considered in Subsection~\ref{RKN}, where $f^{[a]}$ and $f^{[b]}$ are of the form (\ref{eq:f=fa+fbRKN}), the $h$-flow of $F^{[abb]}$ is of the form
\begin{eqnarray*}
  \varphi^{[abb]}_{h}(y,v) = (y, \, v + h^3 g^{[3]}(y)), \qquad
\mbox{where} \quad g^{[3]}(y)=2 g'(y)g(y).
\end{eqnarray*}
This shows in addition that  $\varphi^{[b]}_{h}$ and $\varphi^{[abb]}_{h}$ commute, and that in particular, for arbitrary $b_j,c_j \in \mathbb{R}$,
\begin{eqnarray}
\label{eq:modpotphi}
  \varphi_{b_j\, h}^{[b]} \circ \varphi_{c_j\, h}^{[abb]} \circ \varphi_{b_j\, h}^{[b]}(y,v) = (y, \, v+ 2 h b_j\, g(y) + h^3 c_j\, g'(y)g(y)),
\end{eqnarray}
which is precisely the $h$-flow of the vector field $2 b_j\, F^{[b]} + c_j h^2 F^{[abb]}$. It thus makes sense to construct methods defined as compositions of  $\varphi^{[a]}_{a_j h}$ and maps of the form (\ref{eq:modpotphi}) for $j=1,\ldots,s$.

For Hamiltonian systems $H(p,q)=T(p) + V(q)$ with quadratic kinetic energy $T(p)=\frac12 p^T M p$, 
%can also be written in the form (\ref{rkn.1}). In that case,  $\dot x = f^{[a]}(x)$ (resp.  $\dot x = f^{[b]}(x)$) is the Hamiltonian system with Hamiltonian function $T$ (resp. $V$). 
the vector field
$F^{[abb]}=[[F^{[a]},F^{[b]}],F^{[b]}]$ is the vector field associated to the Hamiltonian function 
$ (\nabla V)^T \nabla V$,  which only depends on the position vector $q$. Thus, (\ref{eq:modpotphi}) is just the $h$-flow of the system with Hamiltonian function
%$\{\{V,T\},V\} = (\nabla V)^T \nabla V$ (indeed, that is a consequence of the correspondence between commutators of two Hamiltonian vector fields and the canonical Poisson bracket $\{H_1,H_2\}$ of the associated Hamitonian functions $H_1$ and $H_2$), which only depends on the position vector $q$. Thus, (\ref{eq:modpotphi}) is just the $h$-flow of the system with Hamiltonian function
\begin{equation}
  \label{eq:modpot}
  2 b_j \, V(q) + c_j\, h^2 (\nabla V(q))^T \nabla V(q),
\end{equation}
which reduces to the potential $V(q)$ of the system for $b_j=1/2$ and $c_j=0$. This explains the term `splitting methods with modified potentials' used in the recent literature~\cite{lopezmarcos97esi,rowlands91ana,wisdom96sco} to refer to splitting methods obtained by composing the $h$-flows of $T$ and modified potentials of the form (\ref{eq:modpot}).

This procedure can be generalized by considering ``modified potentials'' of higher degree in $h$. In particular, the flow $\varphi^{[abbab]}_{h}$ of the vector field
\begin{equation}
F^{[abbab]} =  [[[F^{[a]},F^{[b]}],F^{[b]}],[F^{[a]},F^{[b]}]],  \label{modif.3}
\end{equation}
is of the form $\varphi^{[abbab]}_{h}(y,v) = (y,v+h^5 g^{[5]}(y))$, and similarly
for the vector fields
\begin{eqnarray}
F^{[abbabab]} &=& [[[[F^{[a]},F^{[b]}],F^{[b]}],[F^{[a]},F^{[b]}]],[F^{[a]},F^{[b]}]],  \label{modif.4} \\
F^{[abbaabb]} &=& [[[[F^{[a]},F^{[b]}],F^{[b]}],F^{[a]}],[[F^{[a]},F^{[b]}],F^{[b]}]],  \nonumber
\end{eqnarray}
with $h^5 g^{[5]}(y)$ replaced by $h^7 g^{[7,1]}(y)$ and $h^7 g^{[7,1]}(y)$ respectively. The functions $g^{[5]},g^{[7,1]},g^{[7,2]}$ can be written in terms of $g$ and its partial derivatives (see~\cite{blanes01hor} for more details).

In some applications, the simultaneous evaluation of $g(y)$, $g^{[3]}(y)$, $g^{[5]}(y)$, $g^{[7, 1]}(y)$ and $g^{[7,2]}(y)$ is not substantially more expensive in terms of computational cost than the evaluation of $g(y)$ alone. In that case,
by replacing in the scheme (\ref{eq:splitting}) each $\varphi_{b_i h}^{[b]}$ by the $h$-flow of
\begin{eqnarray}
\label{eq:modpotgen}
  C_h(b_i,c_i,d_i,e_{i1},e_{i2})  & \equiv  &   b_i\, F^{[b]} + h^2 c_i\, F^{[abb]} + h^4 d_i \, F^{[abbab]}  \nonumber \\ 
   & & + h^6 (e_{i,1} F^{[abbabab]} + e_{i,2}\, F^{[abbaabb]}) 
\end{eqnarray}
additional free parameters are introduced to the scheme without increasing too much the computational cost, which allows the construction of more efficient integrators.

Of course, this can be further generalized by considering more general nested commmutators of $F^{[a]}$ and $F^{[b]}$ that gives rise to ``modified potentials''. In that case, higher degree commutators afected by higher powers of $h$ should be added in (\ref{eq:modpotgen}).

Notice that in this case the coefficients $a_i, b_i$ have not to satisfy 
all the order conditions at order $r\geq 3$ and then, the results for negative 
time steps do not apply in this case. As a result, schemes with positive
coefficients do exist. 
In this case, negative coefficients appear in methods of order six \cite{chin05sop}.

\subsection{Methods with processing}
  \label{processing}

Recently, the processing technique has been used
to find composition methods requiring less evaluations than
conventional schemes of order $r$. The idea consists in
enhancing an integrator $\psi_h$ (the \emph{kernel}) with a
parametric map $\pi_h: \mathbb{R}^D \longrightarrow \mathbb{R}^D$
(the \emph{post-processor}) as
\begin{equation}   \label{1.6}
   \hat \psi_h = \pi_h \circ \psi_h \circ \pi_h^{-1}.
\end{equation}
Application of $n$ steps of the new (and hopefully better) integrator
$\hat{\psi}_h$ leads to
\[
  \hat{\psi}_h^n = \pi_h \circ \psi_h^n \circ \pi_h^{-1},
\]
which can be considered as a $h$-dependent change of coordinates in phase space.
Observe that processing is advantageous if $\hat{\psi}_h$
is a more accurate method than $\psi_h$ and, either the cost of $\pi_h$
is negligible or frequent output is not required, since in that case,
it provides the accuracy of $\hat{\psi}_h$
at essentially the cost of the least accurate method $\psi_h$.

The simplest example of a processed integrator is provided in fact
by the St\"ormer--Verlet method. As a consequence of the group
property of the exact flow, we have
\begin{eqnarray}   \label{eq.6.1}
  \mathcal{S}_h^{[2]} & = &  \varphi_{h/2}^{[a]} \circ \varphi_h^{[b]}
  \circ  \varphi_{h/2}^{[a]} = \varphi_{h/2}^{[a]} \circ \varphi_h^{[b]}
       \circ \varphi_{h}^{[a]} \circ \varphi_{-h}^{[a]} \circ
       \varphi_{h/2}^{[a]}
    \nonumber \\
   &  = &  \varphi_{h/2}^{[a]} \circ \chi_{h} \circ
       \varphi_{-h/2}^{[a]}  = \pi_h \circ \chi_{h} \circ   \pi_h^{-1}
\end{eqnarray}
with $\pi_h = \varphi_{h/2}^{[a]}$ and the symplectic Euler method
$\chi_{h} = \varphi_h^{[b]} \circ \varphi_{h}^{[a]}$.
 Hence, applying the basic integrator
$\chi_{h} = \varphi_h^{[b]} \circ \varphi_{h}^{[a]}$ with processing yields a second order of approximation.

Although initially proposed for Runge--Kutta methods \cite{butcher69teo},
the processing technique has proved its usefulness mainly in the
context of geometric numerical integration \cite{hairer06gni}, where
constant step-sizes are widely employed.

We say that the method $\psi_h$ is of \emph{effective order} $r$
if a post-processor $\pi_h$ exists for which $\hat{\psi}_h$ is
of (conventional) order $r$ \cite{butcher69teo}, that is,
\[
  \pi_h \circ \psi_h \circ \pi_h^{-1} = \varphi_h +
    \mathcal{O}(h^{r+1}).
\]
Hence, as the previous example shows, the basic splitting
$\varphi_h^{[b]} \circ \varphi_{h}^{[a]}$ is of effective order 2. Obviously, a method of order
$r$ is also of effective order $r$ (taking $\pi_h = \mathrm{id}$)
or higher, but the converse is not true in general.

The analysis of the order conditions of the method $\hat{\psi}_h$
shows that many of them can be satisfied by $\pi_h$, so that
$\psi_h$ must fulfill a much reduced set of restrictions
\cite{blanes04otn,blanes99siw}. For instance, if the kernel is defined as (\ref{eq:compint}) with a basic first order integrator $\chi_h$ and the post-processor is similarly defined as
\begin{equation}
\label{eq:pp}
  \pi_h =  \chi_{\gamma_{2m} h}\circ
\chi^*_{\gamma_{2s-1}h}\circ
\cdots\circ
\chi_{\gamma_{2}h}\circ
\chi^*_{\gamma_{1}h}
\end{equation}
then, conditions
\begin{equation}
  \label{eq:eoc4}
u_1(\alpha)=1,\quad u_2(\alpha)=u_3(\alpha)=u_4(\alpha)=0
\end{equation}
guarantee that the kernel $\psi_h$ is of effective order four. If in addition the post-processor (\ref{eq:pp}) satisfies 
\begin{equation*}
  u_{1}(\gamma)=0, \quad u_{2}(\gamma)=u_{12}(\alpha), \quad u_{3}(\gamma)=u_{13}(\alpha), \quad u_{12}(\gamma)=u_{112}(\alpha),
\end{equation*}
then the processed integrator (\ref{1.6}) has conventional order
four.  Here, we use the notation $\alpha=(\alpha_1,\ldots,\alpha_{2s})$ and $\gamma=(\gamma_1,\ldots,\gamma_{2m})$ for the coefficients of the kernel and the post-processor respectively.
If in addition the following conditions are fulfilled by the coefficients of the kernel,
\begin{equation*}
  u_{5}(\alpha)=u_{23}(\alpha)=0, \qquad 2 u_{122}(\alpha) + u_{14}(\alpha)+u_{12}(\alpha)^2 = 0,
\end{equation*}
then the kernel $\psi_h$ has at least effective order five. In that
case, the processed method (\ref{1.6}) achieves conventional order
five if in addition, the equalities
\begin{eqnarray*}
   u_{4}(\gamma) \!&=&\! u_{14}(\alpha), \quad u_{13}(\gamma)=u_{113}(\alpha), \\
 u_{112}(\gamma) \!&=&\! u_{1112}(\alpha) + \frac{1}{2} u_{12}(\alpha)^2 - \frac{1}{2} u_{112}(\alpha) - \frac{1}{6} u_{12}(\alpha)
\end{eqnarray*}
hold for the coefficients of the post-processor (\ref{eq:pp}).

Thus, the number and complexity of the conditions to be
verified by the coefficients $\alpha_j$ of a kernel of the form
(\ref{eq:compint}) is notably reduced.  Highly efficient processed
composition methods that take advantage of that have been proposed
\cite{blanes99siw,mclachlan02sm}. Nevertheless, when both the kernel $\psi_h$ and the post-processor $\pi_h$ are constructed as a composition of the form (\ref{eq:compint}) (or (\ref{eq:splitting})), the use of the
resulting processed scheme is not recommended in situations where
intermediate results are required at each step. Indeed, the total number of compositions per step in a processed method (\ref{1.6}) of that form is typically higher than for a non-processed method of comparable accuracy.

To overcome this drawback, in \cite{blanes04otn} a technique has been
developed for obtaining approximations to the post-processor at
virtually cost free and without loss of accuracy. The key idea is to
replace $\pi_h$ by a new map $\tilde{\pi}_h \simeq \pi_h$
obtained from the intermediate stages in the computation of
$\psi_h$. The post-processor $\pi_h$ can safely be replaced
by an approximation $\tilde{\pi}_h$, since the error introduced by the cheap
approximation $\tilde{\pi}_h$ is of a purely local nature~\cite{blanes04otn} (it is not propagated along the evolution, contrarily to the error in
$\pi_h^{-1}$).

In~\cite{blanes04otn}, a general study of the number of independent
effective order order conditions versus the number of conventional
order conditions is presented. In particular, it is shown that, in the
case of kernels of the form (\ref{eq:compint}), the number of
conditions to increase the effective order of the kernel from $k>1$
(resp. $k=1$) to $k+1$ is $n_{k+1}-n_{k}$ (resp. $n_2-n_1+1$), where
each $n_k$ is the cardinal of $L_k$, that is, the number
Lyndon multi-indices of degree $k$. Thus, whereas the total number of
independent conditions to achieve conventional order $r$ is
$n_1+\cdots+n_r$, only $1+n_{r}$ conditions have to be imposed to the
kernel for effective order $r$. If the kernel (\ref{eq:compint}) is
time-symmetric (i.e., if its coefficients satisfy (\ref{eq:sym})),
then there are $N_r=\sum_{i=1}^{q} n_{2j-1}$ independent conditions
for order $r=2q$, and $N_r^* = n_1 + \sum_{i=1}^{q-1}
(n_{2j+1}-n_{2j})$ conditions for effective order $r=2q$. A similar
situation occurs for the total numbers $M_r$ and $M_r^*$ of
conventional and effective order conditions of symmetric kernels of
the form (\ref{eq:compSS}) with (\ref{eq:sym3}) (where the $n_k$ are
replaced by the number $m_k$ in Table~\ref{tabla1}). That
also happens to be true for symmetric kernels of the form
(\ref{eq:splitting}), both in the general case (which is essentially
equivalent to the case of kernels of the form (\ref{eq:compint})) and
in the RKN case considered in Subsection~\ref{RKN}. In
Table~\ref{tb:eoc}, the total number of conditions for conventional
order $r=2q$ for symmetric kernels is compared with the total number
of effective order conditions in three kinds of integrators: (i)
$(N_r,N_r^*)$ for composition (\ref{eq:compint}) of a basic first
order integrator and its adjoint, (ii) $(M_r,M_r^*)$ for compositions
(\ref{eq:compSS}) of a symmetric second order basic integrator,
(iii) $(L_r,L_r^*)$ for splitting integrators in the RKN case.

\begin{table}
\begin{center}
  \begin{tabular}{|c|cccccc|} \hline
  $r$ & 2 & 4 & 6 & 8 & 10 & 12 \\ \hline
$N_r$    & 1 & 3 & 9 & 27& 83 & 269 \\
$N_r^*$  & 1 & 2 & 5 & 14 & 40 & 127 \\ \hline
$M_r$    & 1 &  2 & 4  & 8 & 16 & 33 \\
$M_r^*$  & 1 &  2 & 3  & 5 & 8 & 14 \\ \hline
$L_r$    & 2 &  4 & 8  & 18 & 43 & 112 \\
$L_r^*$  & 2 &  3 & 5  & 10 & 21 & 51 \\ \hline
\end{tabular}
\end{center}
  \caption{Number of conventional and effective
order conditions for symmetric kernels: (i)
$(N_r,N_r^*)$ for composition integrators (\ref{eq:compint}), (ii) $(M_r,M_r^*)$ for compositions (\ref{eq:compSS}) of symmetric second order basic integrator,
(iii) $(L_r,L_r^*)$ for splitting integrators in the RKN case.}
\label{tb:eoc}
\end{table}

\section{A collection of splitting methods}
\label{sec.5}

As we have mentioned before, splitting methods have found application
in many different areas of science during the last decades. It is therefore
not surprising that there is a large number of different schemes available in the
literature. Sometimes, even the same method has been rediscovered several times
in different contexts. Our aim in this section is to offer the reader
a comprehensive overview of the existing methods, by classifying them
into different families and giving the appropriate references where the corresponding
coefficients can be found.

%The simple procedure (\ref{eq:S_h-rec1})-(\ref{eq:S_h-rec2}) to build high order
%methods starting from low order schemes was one of the first breakthroughs in
% a turning point (and probably one of the
%most influential results) in the field of {\it Geometric
%the field of geometric numerical integration. From a low order symmetric method which preserves
%some desired properties of the actual flow, it is possible to build methods preserving
%the same properties at any order and with good error growth properties. Their main
%drawback, however, is the large truncation errors they typically possess, so that their efficiency
%is superior to that of standard schemes only when very long time integrations are
%considered. As we have seen, more efficient schemes can be built by considering other
%types of compositions.

At this point it is important to remark that the efficiency of a method is
measured by taking into account the computational cost required to achieve a
given accuracy (we do not take into account the important
property of the stability of the methods). For instance, given
several methods of order $r$ with different computational cost
(usually measured as the number of stages or evaluations of the functions involved), the most efficient
method does not necessarily correspond to the cheapest method. The
extra cost of some methods can be compensated by an improvement in
the accuracy obtained.

We next present a short review indicating the splitting methods
which have been published in the literature at different orders, with
different number of stages and for several families of problems.

\paragraph{Symmetric compositions of symmetric methods.}

 As we pointed out in
section \ref{s&c}, although by applying recursively the composition
(\ref{eq:S_h-rec1})-(\ref{eq:S_h-rec2}) it is possible to increase the order,
the resulting methods are computationally expensive. To reduce the number
of evaluations the more general composition (\ref{eq:compSS}) may be
considered to achieve a given order $r$. 
 If we choose symmetric compositions ($\beta_{s+1-i}=\beta_i$), then half of the
parameters of the method are fixed, but the order conditions at
even orders are automatically satisfied. In other words, the parameters of a (non-symmetric)
method of order $r=2k$ have to solve a system of
$\sum_{i=1}^{2k}m_i$ equations (see Table~\ref{tb:eoc}), whereas for a
symmetric composition 
$e_r=m_2+m_4+\cdots+m_{2k}$ order conditions are automatically
satisfied if the order conditions at odd orders are fulfilled. In this way,
only $M_r=m_1+m_3+\cdots+m_{2k-1}$ independent order
conditions need to be imposed in the case of symmetric
compositions. Due to this fact, the number of conditions to be
solved  (which is typically the bottleneck in the
numerical search of methods) is reduced considerably when imposing
symmetry. Furthermore, since
$m_{2i-1}<m_{2i}$ then $M_r<e_r$ and symmetric compositions, in
addition to having more favourable geometric properties (due to the
time-symmetric property), usually require smaller number of stages
than their non-symmetric counterparts. 
Taking into account the number $M_r$
(resp. $M_r^*$) of independent conditions to achieve conventional
order $r$ (resp. effective order $r$) from Table~\ref{tb:eoc}, it
is possible to determine the minimum number $s_r=2M_r-1$ of stages
of the integrator (resp. the minimum number $k_r=2M_r^*-1$ of
stages for the kernel) required by a method of order $r$ (resp.
effective order $r$)
%\begin{center}
%\begin{tabular}{|l|c|c|c|c|} \hline
% $r$         & 4  &  6  &  8  &  10 \\
% \hline
% $s_r$  & 3  &  7  & 15 & 31 \\
% \hline
% $k_r$  & 3  & 5  & 9 & 15 \\
% \hline
%\end{tabular}
%\end{center}

 In this way one has to solve a system of $M_r$ or
$M_r^*$ nonlinear polynomial equations with the same number of
unknowns $\beta_i$. The number of real solutions typically increase a good deal with $r$. In general, these
equations have to be solved numerically and getting \emph{all}
solutions is a very challenging task, even for moderate values of
$r$. Once a number of solutions for the parameters $\beta_i$ have
been obtained, there remains to select that solution one expects
will give the best performance when applied on practical problems,
typically by minimizing some objective function. What is the most
appropriate objective function in this case? A frequently used
criterion is to choose the solution which minimizes
$C=\sum_{i=1}^s |\beta_i|$.

If one takes additional stages in (\ref{eq:compSS}), for instance
$s=s_r+2$,  then one has an extra free parameter (notice the
scheme is symmetric and two stages are required to introduce one
parameter). By choosing $\beta_1$ as this free parameter,  then it
is clear that 1-parameter families of solutions are obtained. For
instance, taking $\beta_1=0$ one has the previous solutions and by
continuation it is possible to get several of these  1-parameter
families of solutions, but this procedure does not guarantee to
find all solutions.

Finally, one has to select that solution minimizing the value of
$C$. Of course, additional stages can be introduced and the
process is similar but technically much more involved. This
objective function allows one to find very efficient methods
involving additional stages, although the efficiency of methods
with the same order but different number of stages cannot be
compared from the value obtained for $C$.

When we stop including additional stages in the composition
(\ref{eq:compSS})? Two criteria are possible: (i) when
one has enough stages available to achieve a higher order; (ii)
when the performance of the actual methods constructed with
additional stages do not show a significant improvement in
numerical experiments.

For instance, the simple 4th-order scheme (\ref{suzu1}) can be
improved just by the  5-stage generalized composition
\cite{suzuki90fdo}
\begin{equation}\label{suzukiSS4}
 \psis_{\alpha h}^{[2]}  \circ \psis_{\alpha h}^{[2]}  \circ \psis_{\beta h}^{[2]}
  \circ \psis_{\alpha h}^{[2]} \circ \psis_{\alpha h}^{[2]},
\end{equation}
with $\alpha=1/(4-4^{1/3}), \ \beta=1-4\alpha$, as numerical
experiments clearly indicate. This is a particular case of
(\ref{eq:compSS}) where the value of $C$ reaches a minimum and if we
add two new stages with a new parameter then a 6th-order method
can be obtained.

In Table~\ref{tableSS} we collect some of the most relevant
methods from the literature with different orders and number of
stages. At each order, $r$, we label  the methods by the number of
stages ${\bf s}$  and the reference where this method
can be found. We also include methods obtained by using the
processing technique, which are referred as {\bf P}:${\bf s}$,
where $s$ is the number of stages for the
kernel. We write ${\bf s_1}$-${\bf s_2}$ for indicating that
methods from ${\bf s_1}$ up to ${\bf s_2}$ stages are analyzed in
that particular reference.

\begin{table}
  \centering
\begin{tabular}{|c|c|c|c|} \hline
\multicolumn{4}{|c|}{\bfseries Order}\\ \hline
  % after \\ : \hline or \cline{col1-col2} \cline{col3-col4} ...
4 & 6  & 8 & 10 \\ \hline
  \textbf{3}-\cite{forest90fos,creutz89hoh,yoshida90coh,suzuki90fdo} &  \textbf{7}-\cite{yoshida90coh} &  {\bf 15}-\cite{yoshida90coh,suzuki93hod,mclachlan95otn,kahan97ccf} &  {\bf 31}-\cite{suzuki93hod,kahan97ccf,hairer06gni,sofroniou05dos} \\
   {\bf 5}-\cite{suzuki90fdo,mclachlan95otn} &  {\bf 9}-\cite{mclachlan95otn,kahan97ccf} &  {\bf 17}-\cite{mclachlan95otn,kahan97ccf} &  {\bf 33}-\cite{kahan97ccf,hairer021gni,tsitouras99ato,sofroniou05dos} \\
               & {\bf 11-13}-\cite{sofroniou05dos} &  {\bf 19-21}-\cite{sofroniou05dos} &  {\bf 35}-\cite{hairer021gni,sofroniou05dos} \\
               &  &  {\bf 24}-\cite{calvo93hos} &   \\ \hline \hline

   {\bf P:3-17}-\cite{mclachlan02foh} & {\bf P:5-15}-\cite{blanes06cmf} & {\bf P:9-19}-\cite{blanes06cmf} & {\bf P:15-25}-\cite{blanes06cmf} \\ \hline
\end{tabular}
  \caption{Symmetric compositions of symmetric methods published in the literature. We indicate the
  number of stages (in boldface) and the pertinent reference. Processed methods are preceded by
  \textbf{P}. }\label{tableSS}
\end{table}

\paragraph{Splitting into two parts. Composition of method and its adjoint.}

Next we review methods of the form (\ref{eq:splitting}) (for ODEs
that can be split into two parts) and (\ref{eq:compint}). It is
important to emphasize that, although the order conditions for
both classes of methods are equivalent, the optimization
procedures carried out to identify the most efficient schemes may
differ. In consequence, a particular method optimized for
equations separable into two parts is not necessarily the best choice
for a composition (\ref{eq:compint}), although their performances
are closely related.

Considering, as before, symmetric compositions, i.e.,
$a_{s+1-i}=a_i$,  $b_{s+2-i}=b_i$ in (\ref{eq:splitting}) and
$\alpha_{2s+1-i}=\alpha_i$ in (\ref{eq:compint}), it is easy to
verify, from Table~\ref{tb:eoc}, that the minimum number of stages
required to get a method of order $r$ is $s_r=N_r$ and of
effective order $r$ it is $k_r=N_r^*$.
% is as follows: \
%\begin{center}
%\begin{tabular}{|l|c|c|c|} \hline
% $r$         & 4  &  6  &  8   \\
% \hline
% $s_r$  & 3  &  9  & 27  \\
%  \hline
% $k_r$  & 3  &  5  & 14 \\
%  \hline
%\end{tabular}
%\end{center}

 Note that schemes of order six or higher require more stages
than compositions (\ref{eq:compSS}), and only fourth-order methods
seem promising. Nevertheless, one should recall that by including
additional stages more efficient methods could be obtained. For
instance,  sixth-order methods require at least 9 stages (unless
they are considered as composition of symmetric-symmetric methods
in which case the 9 equations can be solved with only 7 unknowns)
and the coefficients  $a_i,b_i$ or $\alpha_i$ have to solve a
system of eight nonlinear equations (in addition to consistency
conditions). These equations have a very large number of solutions
and it might be the case that one of them could correspond to a
method with very small error terms.

One optimization criterion frequently used when dealing with
composition (\ref{eq:splitting}) is to work with the homogeneous
subspace $\cL_{r+1}=\langle F_{r+1,1},\ldots,F_{r+1,n_{r+1}}\rangle$ (where
by $F_{r+1,i}$ we denote the elements of the basis of the Lie algebra generated by
$F^{[a]}$, $F^{[b]}$ at order $r+1$)
and the leading error term, which can be expressed as
$\sum_{i=1}^{n_{r+1}} c_iF_{r+1,i}$. In this
setting, one selects the solution minimizing
${E}_{r+1}=\left(\sum_{i=1}^{n_{r+1}} |c_i|^2\right)^{1/2}$. This
optimization criterion allows one to compare the performance
of methods with different number of stages by introducing the
effective error, $\mathcal{E}_f =s {E}_{r+1}^{1/r}$, which normalizes with
respect to the number of stages.

For the composition (\ref{eq:compint}), it is not so clear how to
assign a weight to each element of the associated Lie algebra
since their contribution on the error can differ significantly.
One accepted choice consists in minimizing the objective function
$C=\sum_{i=1}^{2s} |\alpha_i|$.

Methods up to order six built by applying this procedure can be
found in the literature. They show for most problems a better
efficiency than compositions (\ref{eq:compSS}) at the same order
when applied to the same class of problems. We collect some of the
most relevant schemes in Table~\ref{tableAB}. As before, we also
include processed methods.

We have not found methods of order eight. In fact, it is an open problem to determine if such
a large system of polynomial equations admits solutions
leading to more efficient methods than those collected in Table~\ref{tableSS}.

\begin{table}
  \centering
\begin{tabular}{|c|c|c|c|} \hline
\multicolumn{4}{|c|}{\bfseries Order}\\ \hline
  % after \\ : \hline or \cline{col1-col2} \cline{col3-col4} ...
  3 & 4 & 6  & 8  \\ \hline
 {\bf 3}-\cite{ruth83aci} & {\bf 3}-\cite{forest90fos,creutz89hoh,yoshida90coh,suzuki90fdo}  &  {\bf 9}-\cite{forest92sol}  &  {\bf 27}-?  \\
              & {\bf 4-5}-\cite{mclachlan95otn} &  {\bf 10}-\cite{blanes02psp} &  \\
             & {\bf 6}-\cite{blanes02psp} &   &   \\ \hline  \hline

  & {\bf P:3,4}-\cite{blanes99siw} & {\bf P:5}-\cite{blanes99siw} & {\bf P:14}-?   \\
  & {\bf P:2-7}-\cite{blanes06cmf} & {\bf P:5-10}-\cite{blanes06cmf} & \\ \hline
\end{tabular}
  \caption{Symmetric composition schemes of the form  (\ref{eq:splitting}) (appropriate when the ODE is
  split in two parts) and
(\ref{eq:compint}) (composition of a method and its adjoint). At order eight, we have not found methods.
They would require at least 27 stages  or
at least 14 stages for processed schemes. The notation is the same as in Table~\ref{tableSS}.}\label{tableAB}
\end{table}

\paragraph{Runge--Kutta--Nystr\"om methods.}

As we have seen in section \ref{RKN}, methods of this class may be
considered as particular examples of composition
(\ref{eq:splitting}). Nevertheless, their wide range of
applicability  to relevant physical problems has originated an
exhaustive search of efficient schemes. Moreover, since in this
case the associated vector fields $F^{[a]}$ and $F^{[b]}$ have
different qualitative properties, methods with different features
may be found in the literature. Thus, one may find
non-symmetric methods of the form
\begin{eqnarray}  \label{ab-ba}
 \mbox{ \textit{AB} } \equiv \, 
 \psi_h & = & \varphi^{[a]}_{a_{s}h}\circ \varphi^{[b]}_{b_{s} h}\circ
 \cdots\circ
 \varphi^{[a]}_{a_{1}h}\circ \varphi^{[b]}_{b_{1}h} \nonumber \\
 \mbox{ \textit{BA} } \equiv \, 
 \psi_h & = & \varphi^{[b]}_{b_{s} h} \circ \varphi^{[a]}_{a_{s}h}\circ
 \cdots\circ
 \varphi^{[b]}_{b_{1}h}\circ \varphi^{[a]}_{a_{1}h}
\end{eqnarray}
where $AB$ and $BA$ are conjugate to each other, leading to the same
performance. However, to take profit of the FSAL (First Same As
Last) property, we can consider the following non equivalent
compositions
\begin{equation}  \label{aba}
 \mbox{ \textit{ABA} } \equiv
 \psi_h = \varphi^{[a]}_{a_{s+1}h}\circ \varphi^{[b]}_{b_{s} h}\circ
 \varphi^{[a]}_{a_{s}h}\circ
 \cdots\circ
 \varphi^{[b]}_{b_{1}h}\circ
 \varphi^{[a]}_{a_{1}h}
\end{equation}
and
\begin{equation}
  \label{bab}
  \mbox{ \textit{BAB} } \equiv
\psi_h = \varphi^{[b]}_{b_{s+1} h}\circ
 \varphi^{[a]}_{a_{s}h}\circ \varphi^{[b]}_{b_{s} h}\circ
 \cdots\circ
 \varphi^{[a]}_{a_{1}h} \circ \varphi^{[b]}_{b_{1}h}.
\end{equation}
The symmetric case ($a_{s+2-i}=a_i,b_{s+1-i}=b_i$ for the
composition $ABA$ and $b_{s+2-i}=b_i,a_{s+1-i}=a_i$ for the
composition $BAB$) has also been proposed in this setting to get more efficient schemes.
In this case, the minimum number of stages is, from Table~\ref{tb:eoc}, $s_r=L_r-1$ 
(resp. $k_r=L_r^*-1$), to get a method of order $r$ (resp. effective order $r$).
For non-symmetric compositions this minimum number can be obtained
from Table~\ref{tb:number-oc-RKN}.
%is as follows: \
%\begin{center}
%\begin{tabular}{|l|c|c|c|c|}  \hline
% $r$         & 4  &  5 &  6  &  8   \\
% \hline
% $s_r$  & 3  &  5  &  7  & 17  \\
% \hline
% $k_r$  & 3  &  3 &  4  & 9  \\
% \hline
%\end{tabular}
%\end{center}

 Highly efficient methods up to order six have been published. In
Table~\ref{tableRKN} we collect the most representative within
this class. We add {\bf S} or {\bf N} to distinguish symmetric from non-symmetric schemes
and the subindex \textit{AB}, \textit{ABA} and \textit{BAB} to denote compositions
(\ref{ab-ba}), (\ref{aba}) and (\ref{bab}), respectively.
Processed methods have also been included.

To achieve order eight, the coefficients $a_i$, $b_i$ in a
non-processed scheme have to solve (in addition to consistency) a
system of 16 nonlinear equations. A large number of solutions
could exist, although, as far as we know, only  one attempt to
solve these equations has been reported \cite{okunbor94eoe} (the
performance of such method was not clearly superior
symmetric-symmetric methods).

\begin{table}
  \centering
\begin{tabular}{|c|c|c|c|} \hline
\multicolumn{4}{|c|}{\bfseries Order}\\ \hline
  % after \\ : \hline or \cline{col1-col2} \cline{col3-col4} ...
 4 & 5 & 6  & 8  \\ \hline
 {\bf 3S}-\cite{forest90fos,creutz89hoh,yoshida90coh,suzuki90fdo} & {\bf 5N}$_{ABA}$-\cite{okunbor94crk}  &  {\bf 7S}$_{ABA}$-\cite{forest92sol,okunbor94crk} & {\bf 17S}$_{ABA}$-\cite{okunbor94eoe} \\
 {\bf 4N}$_{AB}$-\cite{mclachlan92tao}      & {\bf 6N}$_{AB}$-\cite{mclachlan92tao} &  {\bf 7S}$_{BAB}$-\cite{forest92sol} &  \\
 {\bf 4N}$_{BAB}$-\cite{calvo93tdo}     & {\bf 6N}$_{AB}$-\cite{chou00ofs} &     {\bf 8-15S}$_{ABA,BAB}$-\cite{blanes02psp}    &  \\
 {\bf 4-5S}$_{ABA}$-\cite{mclachlan95otn}  &         &  {\bf 7,11S}$_{BAB}$-\cite{blanes02psp}  &  \\
 {\bf 5S}$_{BAB}$-\cite{blanes01smf}     &         &  &  \\
 {\bf 6S}$_{ABA,BAB}$-\cite{blanes02psp}     &         &   &  \\ \hline\hline
 {\bf P:2N}$_{AB}$-\cite{blanes99siw} &   & {\bf P:4-6S}$_{ABA,BAB}$-\cite{blanes01hor} & {\bf P:9S}$_{ABA}$-\cite{blanes01hor}  \\
  &     & {\bf P:7S}$_{BAB}$-\cite{blanes01nfo} & {\bf P:11S}$_{BAB}$-\cite{blanes01nfo}  \\ \hline
\end{tabular}
  \caption{RKN splitting integrators. Since the role of the flows
$\varphi_t^{[a]}$ and $\varphi_t^{[b]}$ is not interchangeable
here, we distinguish symmetric {\bf S} and non-symmetric {\bf N}
compositions with a subindex $AB,ABA,BAB$ for the compositions
(\ref{ab-ba}), (\ref{aba}) and (\ref{bab}). As usual, processed
methods are preceded by {\bf P}.}\label{tableRKN}
\end{table}

In \cite{blanes02psp} the authors have carried out a detailed
analysis of the order conditions for symmetric compositions $ABA$
and $BAB$. In this work 4th-order methods from 3 to 6 stages, and
also 6th-order methods from 7 to 14 stages are analysed. The
integrators selected perform extraordinarily well indeed. For
instance, on the H\'enon--Heiles Hamiltonian (\ref{HenonHeiles})
the 4th-order 6-stage method is more accurate (at constant work)
than leapfrog in a wide range of step sizes, whereas its global
error is about $0.00175$ times that of the classical 4th-order
Runge--Kutta method. In consequence, its computational cost for a
given error is about $0.31$. This has to be compared with the
composition (\ref{suzu1}) based on leapfrog, which have truncation
errors about 10 times larger than the classical Runge--Kutta
scheme.

On the other hand, as we have seen in subsection \ref{modified},
the particular structure of problem (\ref{rkn.1}) allows one to
use modified potentials in compositions (\ref{ab-ba})-(\ref{bab}).
This is appealing when the evaluation of such modified potentials
is not particularly costly.
%As mentioned there, the associated
%vector fields $F^{[a]}$ and $F^{[b]}$ verify
%$[F^{[b]},[F^{[b]},[F^{[b]},F^{[a]}]]] =  0$, so that $F^{[b]}$
%and $[F^{[b]},[F^{[b]},F^{[a]}]]$ commute and, frequently, both
%can be computed exactly and efficiently.
In such circumstances one may
replace in (\ref{ab-ba})-(\ref{bab}) flows associated to
$hb_iF^{[b]}$ with the corresponding to $h
C_h(b_i,c_i,d_i,e_{i1},e_{i2})$, as given in (\ref{eq:modpotgen}). The
coefficients $c_i$, $d_i$, etc. can be used to solve some order
conditions, so that methods with a reduced number of stages can be
obtained. We emphasize that these schemes are of interest when the
extra cost due to the modified potentials is moderate, as is the
case in many problems arising in classical and quantum mechanics.
In Table~\ref{tableRKNm} we collect some relevant methods we have
found in the literature, both processed and non-processed.

\begin{table}
  \centering
\begin{tabular}{|c|c|c|c|} \hline
\multicolumn{4}{|c|}{\bfseries Order}\\ \hline
  % after \\ : \hline or \cline{col1-col2} \cline{col3-col4} ...
 3 & 4 & 6  & 8  \\ \hline
 {\bf 2N}$_{AB}$-\cite{ruth83aci} & {\bf 2S}$_{ABA,BAB}$-\cite{koseleff94fcf,chin97sif} & {\bf 4,5S}$_{ABA,BAB}$-\cite{omelyan02otc} & {\bf 11S}$_{ABA,BAB}$-\cite{omelyan02otc} \\
                      & {\bf 4S}$_{ABA,BAB}$-\cite{suzuki95hep} &   &  \\
                      & {\bf 3,4S}$_{ABA,BAB}$-\cite{chin01fog,omelyan02otc} &   &  \\
 \hline\hline
    & {\bf P:1S}$_{BAB}$-\cite{takahashi84mcc,rowlands91ana,wisdom96sco,blanes99siw}  &  {\bf P:3S}$_{ABA,BAB}$-\cite{blanes99siw} & {\bf P:4S}$_{ABA}$-\cite{blanes01hor}  \\
                   & {\bf P:2S}$_{BAB}$-\cite{lopezmarcos97esi}   &  & {\bf P:5S}$_{BAB}$-\cite{blanes01hor, blanes01nfo}  \\ \hline
\end{tabular}
  \caption{RKN splitting methods with modified potentials. Schemes are coded as in Table
  \ref{tableRKN}.}\label{tableRKNm}
\end{table}

\paragraph{Methods for near-integrable systems.}

As we have seen in section \ref{Q&I}, splitting methods designed
for equation (\ref{qi.1}) have typically two relevant parameters:
$h$ (the step size) and $\varepsilon$ (the size of the
perturbation). In consequence, the dominant error in a given
scheme depends on their relative size, and this depends usually on
the particular problem considered (and sometimes even on the
initial conditions). For this reason, a number of methods at
different orders in both parameters $h$ and $\varepsilon$ are
found in the literature. We collect some of them in
Table~\ref{tableNI}. Here the notation is a bit clumsy: a method of order (n,4), say,
means that the exact and the modified vector fields, i.e.,
$hF$ and $\log(\Psi(h))$ in (\ref{c.1}),
differ in
terms $\mathcal{O}(\varepsilon h^{n+1} + \varepsilon^2 h^5 +
\cdots)$, whereas in a method (7,6,4) this difference is
$\mathcal{O}(\varepsilon h^{8} + \varepsilon^2 h^7 + \varepsilon^3
h^5 + \cdots)$. In both cases, the order of consistency in the
limit $h \rightarrow 0$ is four, but the last method incorporates
more terms in the asymptotic expansion of the error.

In \cite{mclachlan95cmi} both families ($ABA$ and $BAB$) of
symmetric $(2s,2)$ schemes for $s\leq 5$ with positive
coefficients are proposed which are about three times more
accurate (at constant work) than leapfrog, whereas in
\cite{laskar01hos} a systematic study of $(2s,2)$ methods is
carried out, obtaining new schemes up to $s=10$ with positive
coefficients.

In some near-integrable problems, the identity
$[[[F^{[b]},F^{[a]}],F^{[b]}],F^{[b]}] =  0$
 still holds, where
$F^{[i]}$ is the vector field associated to $f^{[i]}$, $i=a,b$, in
(\ref{qi.1}). This takes place, in particular, in Hamiltonian
problems $H = H_0 + \varepsilon H_1$ where $H_0$ is quadratic in
the kinetic energy and $\varepsilon H_1$ depends only on the
coordinates (e.g. examples (\ref{HenonHeiles})-(\ref{wave1}) can
be split in this way, where $H_0$ is the harmonic oscillator or
the Kepler problem and $H_1$ depends only on the coordinates). In
consequence, the previous techniques used to obtain RKN methods
still apply here, as well as the inclusion of flows of modified
potentials in the composition.

In Table \ref{tableNI} we separate, as usual, non-processed from
processed schemes (preceded by {\bf P}). In the later case we also
include methods applicable when
$[[[F^{[b]},F^{[a]}],F^{[b]}],F^{[b]}] =  0$
 (second row) and schemes
with modified potentials (last two rows of processed methods).

\begin{table}
  \centering
\begin{tabular}{|c|c|c|} \hline
\multicolumn{3}{|c|}{\bfseries Order}\\ \hline
  % after \\ : \hline or \cline{col1-col2} \cline{col3-col4} ...
 (n,2) & (n,4) & (n,5)    \\ \hline
{\bf 1}(2,2){\bf S}-\cite{wisdom91smf} & {\bf 4}(6,4){\bf N}$_{BAB}$-\cite{mclachlan95cmi}  &    \\
 {\bf n}(2n,2){\bf N}$_{ABA,BAB}$-\cite{mclachlan95cmi,laskar01hos} & {\bf 5}(8,4){\bf S}$_{ABA,BAB}$-\cite{mclachlan95cmi} & \\ \hline\hline
 {\bf P:1}(32,2)-\cite{wisdom96sco} &  {\bf P:3}(7,6,4){\bf S}$_{ABA}$-\cite{blanes00psm} &  \\ \hline
  & {\bf P:2}(6,4){\bf S}$_{AB}$-\cite{blanes00psm} & {\bf P:3}(7,6,5){\bf S}$_{AB}$-\cite{blanes00psm}  \\ \hline
  & {\bf P:1}(6,4){\bf S}$_{ABA}$-\cite{blanes00psm} & {\bf P:2}(7,6,5){\bf S}$_{AB}$-\cite{blanes00psm}  \\
  & {\bf P:n}(n,4){\bf S}$_{ABA}$-\cite{laskar01hos} &   \\ \hline
\end{tabular}
  \caption{Splitting methods for near-integrable systems.
For processed methods we also include methods applicable when
$[[[F^{[b]},F^{[a]}],F^{[b]}],F^{[b]}] =  0$  (second row) and schemes
with modified flows (last two rows of processed
methods).}\label{tableNI}
\end{table}

\section{Preserving properties and backward error analysis}
\label{bea}

Much insight into the long-time behavior of splitting methods (including preservation of invariants and
structures in the phase space) can be gained by applying backward error analysis techniques. We
will summarize here some of the main issues involved and refer the reader to \cite{hairer06gni} for a
detailed treatment of the theory.

When we analyzed in the Introduction the symplectic Euler scheme as applied to the simple harmonic
oscillator,
we associated its good qualitative properties with the fact that
the numerical solution can be interpreted as the exact solution of a perturbed Hamiltonian system. This
remarkable feature constitutes a simple illustration of the insight provided by
backward error analysis (BEA) in this setting. More generally, suppose that we apply the splitting
method (\ref{eq:splitting}) to solve equation (\ref{harmonic2}). Then the corresponding
numerical solution at time $t=h$ is given by
\[
    x(h) = K(h) x_0 \equiv \e^{b_{s+1} h B} \, \e^{a_s h A} \, \e^{b_s h B} \, \cdots
\e^{b_2 h B} \, \e^{a_1 h A} \, \e^{b_{1} h B} x_0,
\]
where the so-called stability matrix $K(h)$  is given explicitly by
\[
   K(h) = \left( \begin{array}{cr}
                          1  &  0 \\
                         -b_{s+1} h & 1
                    \end{array} \right) \,
                \left( \begin{array}{cr}
                          1  &  a_s h \\
                          0 & 1
                    \end{array} \right) \,  \cdots
                \left( \begin{array}{cr}
                          1  &  a_1 h \\
                          0 & 1
                    \end{array} \right) \,
               \left( \begin{array}{cr}
                          1  &  0 \\
                         -b_{1} h & 1
                    \end{array} \right).
\]
In this way, one gets
\[
    K(h) =      \left( \begin{array}{cr}
                          K_1(h)  &  K_2(h) \\
                          K_3(h) &   K_4(h)
                    \end{array} \right)
\]
where $K_1(h)$, $K_4(h)$ (respectively, $K_2(h)$, $K_3(h)$) are even (resp. odd) functions and
$\det K(h) = 1$. As a matter of fact, any splitting method is uniquely determined by its stability matrix, so that
the analysis can be carried out only with $K(h)$ \cite{blanes07otl}. 
If in addition the splitting method is symmetric then
$K(h)^{-1} = K(-h)$ and we can write
\[
    K(h) =      \left( \begin{array}{cr}
                          p(h)  &  K_2(h) \\
                          K_3(h) &   p(h)
                    \end{array} \right)
\]
where $p(h) = \frac{1}{2} \mbox{tr} (K(h)) = \frac{1}{2} (K_1(h) + K_4(h))$. It can be shown that
the matrix $K(h)$ is stable for a given
$h \in \mathbb{R}$, i.e., $K(h)^n$ is bounded for all the iterations $n$, if and only if there exist real functions
$\phi(h), \gamma(h)$ such that $p(h) = \cos(\phi(h))$ and $K_2(h) = -\gamma(h)^2 K_3(h)$. In that
case
\[
  K(h) = \left(  \begin{array}{cc}
       \cos(\phi(h)) &  \gamma(h) \sin(\phi(h)) \\
       -\frac{\sin(\phi(h))}{\gamma(h)} &  \cos(\phi(h))
         \end{array}  \right) = \exp \left(  \begin{array}{cc}
                                                                    0  &  \gamma(h) \phi(h)  \\
                                                                    -\frac{\phi(h)}{\gamma(h)}  &  0
                                                                 \end{array}  \right)
\]
where, by consistency, $\phi'(0) = 1$ and $\gamma(0) = 1$, whereas symmetry imposes 
$\phi(-h) = -\phi(h)$ and $\gamma(-h) = \gamma(h)$.

This result implies, in particular, that the numerical solution $(q_n, p_n)$ at time $t_n = n h$ obtained
by applying the splitting method to the linear system (\ref{harmonic2}) verifies
\[
  \left(  \begin{array}{c}
         q_n  \\
         p_n
         \end{array} \right) = \left(  \begin{array}{cc}
                                          \cos(t_n \tilde{\omega})  &  \gamma(h) \sin(t_n \tilde{\omega}) \\
                                 -\gamma(h)^{-1} \sin(t_n \tilde{\omega})  &  \cos(t_n \tilde{\omega})
                                 \end{array}  \right) \,
       \left(  \begin{array}{c}
         q_0  \\
         p_0
         \end{array} \right)
\]
for values of $h$ such that $K(h)$ is stable. Here $\tilde{\omega} = \phi(h)/h$. Equivalently,
\[
  q_n = \tilde{q}(t_n), \qquad p_n = (\phi(h) \gamma(h)/h)^{-1} \frac{d}{dt} \tilde{q}(t_n),
\]
where $\tilde{q}(t)$ is the exact solution of
\[
 \frac{d^2 }{dt^2} \tilde{q} + \tilde{\omega}^2 \tilde{q} = 0
\]
with initial condition $\tilde{q}(0) = q_0$, $\tilde{q}'(0) =  (\phi(h) \gamma(h)/h) p_0$. In other words,
the numerical solution provided by the splitting method is the exact solution of a
harmonic oscillator  with frequency $\tilde{\omega} \approx 1$, i.e., of a system of equations satisfying the
same geometric properties as the original system. The existence of such a backward error interpretation
has direct implications for the qualitative behavior of the numerical solution, as well as for its global error.

The main idea can be extended to an arbitrary non-linear ODE
(\ref{eq.1.1}).  Recall from Subsection~\ref{ss:diffop} that each
integrator $\psi_h$ has associated a series $\Psi(h)=I+h \Psi_1 + h^2 \Psi_2
+\cdots$ of differential operators acting on smooth functions $g \in
\cinf$, and its formal logarithm $\log(\Psi(h))$ is a series of
vector fields (viewed as first order differential operators)
$\log(\Psi(h)) = h F_1 + h^2 F_2 + \cdots$. For $g\in \cinf$, the
result of acting each $F_k$ on $g$ is of the form $F_k[g] = g'(x)
f_k(x)$, for a certain smooth map $f_k: \mathbb{R}^D \longrightarrow
\mathbb{R}^D$. Now, consider the
\emph{modified differential equation} (defined as a formal series in powers of $h$)
\begin{equation}   \label{bea.1}
   \tilde{x}' = f_h(\tilde{x}) \equiv f(\tilde{x}) + h f_2(\tilde{x}) + h^2 f_3(\tilde{x}) + \cdots
\end{equation}
associated to the integrator $\psi_h$. Then one has that
$x_n = \tilde{x}(t_n)$, with $t_n = n h$, which allows studying the long-time behaviour of the numerical integrator by analysing the solutions of the system (\ref{bea.1}) viewed as a small perturbation of the original system (\ref{eq.1.1}).
This allows one to get important qualitative information about the numerical solution. In particular,
\begin{itemize}
 \item for symmetric methods, the modified differential equation only contains even powers of $h$;
 \item for volume-preserving methods applied to a divergence-free dynamical system, the
 modified equation is also divergence-free;
 \item for symplectic methods applied to a Hamiltonian system, the modified differential equation
 is (locally) Hamiltonian.
\end{itemize}
In the particular case of a symplectic integration method, this means that there exist smooth
functions $H_j: \mathbb{R}^{2d} \longrightarrow \mathbb{R}$ for $j=2, 3, \ldots$, such that
$f_j(x) = J \nabla H_j(x)$, where $J$ is the canonical symplectic matrix. In consequence,
there exists a modified Hamiltonian of the form
\begin{equation}   \label{bea.2}
   \tilde{H}(q,p) = H(q,p) + h H_2(q,p) + h^2 H_3(q,p) + h^3 H_4(q,p) + \cdots
\end{equation}
such that the modified differential equation is given by
\[
    q' = \nabla_p \tilde{H}(q,p), \qquad p' = -\nabla_q \tilde{H}(q,p).
\]
Of course, if the method has order $r$, say, then  $H_i = 0$ for $i \le r$ in (\ref{bea.2}). In
other words, the modified Hamiltonian has the form $\tilde{H} = H +  h^r H_{r+1} + \cdots$.
In particular, for the St\"ormer-Verlet method (\ref{leapfrog}) applied to the Hamiltonian
$H(q,p) = T(p) + V(q)$, one has
\[
   \tilde{H} = H + h^2 \left( -\frac{1}{24} V_{qq}(T_p,T_p) + \frac{1}{12} T_{pp}(V_q,V_q) \right)
      + \cdots
\]

Apart from the linear case analyzed before, the series in (\ref{bea.1}) does not converge in general.
To make this formalism rigorous, one has to give bounds on the coefficient functions $f_j(x)$ of
the modified equation so as to determine an optimal truncation index and finally one has to
estimate the difference between the numerical solution $x_n$ and the exact solution $\tilde{x}(h)$
of the modified equation.

These estimates constitute in fact the basis for rigorous statements
about the long term behavior of the numerical solution.  For instance, this theory allows one to
proof rigorously that a symplectic numerical method of order $r$ with constant step size $h$
applied to a Hamiltonian system  $H$ verifies that $H(x_n) = H(x_0) + \mathcal{O}(h^r)$ for
exponentially long time intervals \cite{hairer06gni}.

On the other hand, since the modified differential equation of a numerical scheme depends
explicitly on the step size used, one has a different modified equation each time the step size $h$
is changed. This fact seems to be the reason of the poor long time behavior observed in practice when
a symplectic scheme is implemented directly with a standard variable step-size strategy.

\section{Special methods for special problems}

\subsection{Splitting methods for linear systems}

Suppose one is interested in solving numerically the differential equations arising from the
generalized harmonic oscillator with Hamiltonian function (\ref{harmonicGeneral}).
Although RKN methods with modified potentials can be always used for this purpose,
we will see in the sequel that
the particular structure of this system allows one to design specially
tailored schemes which are orders of magnitude more efficient than other integrators frequently used
in the literature.

At this point, the reader could reasonably ask about the convenience of designing new numerical methods for the
harmonic oscillator (\ref{harmonicGeneral}). It turns out, however, that efficient splitting methods for this system
can be of great interest for the numerical treatment of partial differential equations appearing
in quantum mechanics, optics and electrodynamics previously discretized in space.

Suppose, in particular, that we have to solve numerically the time dependent Schr\"odinger equation
(\ref{Schrodinger})
with initial wave function $\psi(x,0)=\psi_0(x)$. We can write
(\ref{Schrodinger}) as
\begin{equation}  \label{Schr2}
 i \frac{\partial}{\partial t} \psi = \left( T(P) + V(X) \right) \psi,
\end{equation}
where $\displaystyle T(P)=\frac{1}{2 m} P^2$, and the operators
$X, \, P$ are defined by their actions on $\psi (x,t)$ as
\[
  X\psi (x,t) = x\psi (x,t), \qquad\quad
  P\psi (x,t) = -i \nabla \psi (x,t).
\]
For simplicity, let us consider the one-dimensional
problem and suppose that it is defined in a given interval $ \, x
\in [x_0, x_N] $ ($\psi(x_0,t)=\psi(x_N,t)=0$ or it has periodic
boundary conditions). A common procedure consists in taking first
a discrete spatial representation of the wave function
$\psi(x,t)$: the interval is split in $N$ parts of length $\Delta
x = (x_{N}-x_{0})/N$ and the vector $\mathbf{u} = (u_{0}, \ldots,
u_{N-1})^T \in \mathbb{C}^N$ is formed, with $u_{n} =
\psi(x_{n},t)$ and $x_{n} = x_{0} + n \Delta x$,
$n=0,1,\ldots,N-1$. The partial differential equation
(\ref{Schr2}) is then replaced by the $N$-dimensional linear ODE
\begin{equation}   \label{lin1}
  i \frac{d }{dt} \mathbf{u}(t) = \mathbf{H} \, \mathbf{u}(t),  \qquad
    \mathbf{u}(0) = \mathbf{u}_{0} \in \mathbb{C}^N,
\end{equation}
where $\mathbf{H} \in \mathbb{R}^{N\times N}$ represents the (in
general Hermitian) matrix associated with the Hamiltonian
\cite{feit82sot}. The formal solution of equation (\ref{lin1}) is given
by $\mathbf{u}(t) = \e^{-i t \mathbf{H}} \mathbf{u}_{0}$, but to
exponentiate this $N \times N$ complex and full matrix can be
prohibitively expensive for large values of $N$, so in practice
other methods are preferred.

In general $\mathbf{H}=\mathbf{T}+\mathbf{V}$, where $\mathbf{V}$
is a diagonal matrix associated with the potential energy $V$ and
$\mathbf{T}$ is a full matrix related to the kinetic energy $T$.
Their action on the wave function vector is obtained as follows.
The potential operator being local in this representation, one has
$(\mathbf{V u})_{n} = V(x_{n}) u_{n}$ and thus the product
$\mathbf{V  u}$ requires to compute $N$ complex multiplications.
Since periodic boundary conditions are assumed, for the kinetic
energy one has ${\bf T}\, {\bf u} = \mathcal{F}^{-1} {\bf D}_T
\mathcal{F}{\bf u}$, where $\mathcal{F}$ and $\mathcal{F}^{-1}$
correspond to the forward and backward discrete Fourier transform,
and ${\bf D}_T$ is local in the momentum representation (i.e., it
is a diagonal matrix). The transformation $\mathcal{F}$ from the
discrete coordinate representation to the discrete momentum
representation (and back) is done via the fast Fourier transform
(FFT) algorithm, requiring $\mathcal{O}(N \log N)$ operations. It is therefore possible
to use the methods of subsection \ref{s&c} with this splitting.

There are other ways, however, of using splitting techniques in this context. To this end,
notice that $\e^{-it{\bf H}}$ is not only unitary, but also symplectic with canonical coordinates
$\mathbf{q} = \mbox{Re}(\mathbf{u})$ and momenta $\mathbf{p} = \mbox{Im}(\mathbf{u})$. Thus,
equation (\ref{lin1}) is equivalent to
\cite{gray96sit,gray94chs}
\begin{equation}   \label{clas1}
  \frac{d }{dt} \mathbf{q} = \mathbf{H} \, \mathbf{p}, \, \qquad  \qquad
  \frac{d }{dt} \mathbf{p} = - \mathbf{H} \, \mathbf{q},
\end{equation}
where $\mathbf{H \, q}$ and $\mathbf{H \, p}$ require both a
real-complex FFT and its inverse. In addition, system
(\ref{clas1}) can be seen as the classical evolution equations
corresponding to the Hamiltonian function (\ref{harmonicGeneral}) with $M = N = \mathbf{H}$.
Thus, efficient schemes for solving numerically the generalized harmonic oscillator can be applied
directly to this problem. 
Also the Maxwell equations (\ref{Maxwell})  in an isotropic, lossless and source free medium, when
they are previously discretized in space have a similar structure \cite{rieben04hos}. 
In consequence, numerical methods
of this class are well adapted for their numerical treatment.

Clearly, one may write
\begin{equation}  \label{clasic}
 \frac{d}{dt} \left\{
 \begin{array}{c}
  \mathbf{q} \\
  \mathbf{p}  \end{array} \right\} =  \left(
 \begin{array}{ccc}
  \mathbf{0}   & \, & \mathbf{H}  \\
  -\mathbf{H} & \, &  \mathbf{0}  \end{array} \right)   \left\{
 \begin{array}{c}
  \mathbf{q} \\
  \mathbf{p}  \end{array} \right\} = ( \mathbf{A}+ \mathbf{B}) \left\{
 \begin{array}{c}
  \mathbf{q} \\
  \mathbf{p}  \end{array} \right\},
\end{equation}
with the $2N \times 2N$ matrices $\mathbf{A}$ and $\mathbf{B}$
given by
\begin{equation*}    \label{eq.4a}
  \mathbf{A} \equiv \left(  \begin{array}{ccc}
               \mathbf{0}   & \, &  \mathbf{H}  \\
              \mathbf{0}    & \, &  \mathbf{0}  \end{array} \right),
              \qquad\qquad
  \mathbf{B} \equiv \left(  \begin{array}{ccc}
               \mathbf{0}   & \, & \mathbf{0}  \\
              -\mathbf{H}    & \, & \mathbf{0}  \end{array} \right).
\end{equation*}
The evolution operator corresponding to (\ref{clasic}) is
\begin{equation}  \label{exact-Sch}
   \mathbf{O}(t) = \left(
 \begin{array}{rcr}
   \cos(t \mathbf{H})   &  & \sin(t \mathbf{H})  \\
  -\sin(t \mathbf{H})  &   & \cos(t \mathbf{H})  \end{array}
  \right),
\end{equation}
which is an orthogonal and symplectic $2N \times 2N$ matrix. As
before, its evaluation is computationally very expensive and thus
some approximation is required. The usual procedure is to split
the whole time interval into $M$ steps of length $h=t/M$, so that
$\mathbf{O}(t)=[\mathbf{O}(h)]^M$, and then approximate
$\mathbf{O}(h)$ acting on the initial condition at each step.

In this respect, observe that
\begin{equation*} \label{eA_eB}
  \e^{\mathbf{A}} = \left(
 \begin{array}{ccc}
  \mathbf{I}   & \,  & \mathbf{H}  \\
  \mathbf{0}  & \,  &  \mathbf{I}   \end{array} \right), \qquad \qquad
   \e^{\mathbf{B}} = \left(
 \begin{array}{ccc}
   \mathbf{I}   & \, & \mathbf{0} \\
  -\mathbf{H}  &  \, & \mathbf{I}  \end{array} \right)
\end{equation*}
and the cost of evaluating the action of $\e^{\mathbf{A}}$ and
$\e^{\mathbf{B}}$ on $\mathbf{z} = (\mathbf{q},\mathbf{p})^T$ is
essentially the cost of computing the products $\mathbf{H \, p}$
and $\mathbf{H \, q}$, respectively. It makes sense, then, to use
splitting  methods of the form (\ref{eq:splitting}), which in this context read
\begin{equation} \label{compos}
%  e^{\tau ({\bf A+B})}
 {\bf O}_n(h) =
% \prod_{i=1}^k  e^{h a_i{\bf A}} \, e^{h b_i{\bf B}}
  \e^{h b_{s+1}{\bf B}} \, \e^{h a_s{\bf A}} \ \cdots \
  \e^{h b_2{\bf B}} \, \e^{h a_1{\bf A}} \, \e^{h b_1{\bf B}}.
\end{equation}

Several methods with different orders have been constructed along
these lines indeed \cite{gray96sit,liu05oos,zhu96nmw}. Of particular relevance
are the schemes presented in \cite{gray96sit}, since only $s=r$
exponentials $\e^{h a_{i} \mathbf{A}}$ and $\e^{h b_{i}
\mathbf{B}}$ are used to achieve order $r$ for $r=4,6,8,10$ and
$12$. By contrast, in a general composition (\ref{compos}) the
minimum number $s$ of exponentials $\e^{h a_{i} \mathbf{A}}$ and
$\e^{h b_{i} \mathbf{B}}$ (or stages) required to attain order $\
r=8,10 \ $ is $k=15,31$, respectively \cite{hairer06gni,kahan97ccf}.

Furthermore, one can use processing to reduce even more the number of exponentials.
A different approach can also be followed, however: to take a number of stages
larger than strictly necessary to
solve all the order conditions to improve the efficiency and
stability of the resulting schemes. The idea is to use the extra cost to reduce the size of the error
terms, enlarge the stability interval and achieve therefore a higher efficiency but without raising
the order.

Kernels with up to $19$,
$32$ and $38$ stages have been proposed, and
for each kernel  the corresponding coefficients $a_i,b_i$ have been determined
according to two different criteria. The first set of solutions is
taken so as to provide methods of order $r=10$, $16$ and $20$. The
second set of coefficients bring highly accurate \emph{second}
order methods with an enlarged domain of stability. A more detailed treatment can be
found in \cite{blanes06sso,blanes07otl}.

\subsection{Splitting methods for non-autonomous systems}

So far we have considered the problem of designing splitting methods for the numerical
integration of autonomous differential equations (\ref{eq.1.1}). As we have shown, there
are a large number of schemes of different orders in the literature, and some of them are
particularly efficient when the system possesses some additional structure, e.g., for the
second-order differential equation $y'' = g(y)$ and the generalized harmonic oscillator
(\ref{harmonicGeneral}). In this section we will review two different strategies to apply the
splitting schemes when
there is an explicit time dependency in the original problem.

To fix ideas, let us assume that our system is non-autonomous and can be split as
\begin{equation}    \label{naut.1}
   x^{\prime} = f(x,t) = f^{[a]}(x,t) + f^{[b]}(x,t),
   \qquad x(0) = x_0.
\end{equation}
The first, most obvious procedure consists in taking $t$ as a new coordinate, so that (\ref{naut.1})
is transformed into an equivalent autonomous equation to which standard splitting algorithms can
be applied. More specifically, equation (\ref{naut.1}) is equivalent to the
enlarged system
\begin{equation}  \label{naut.2}
    \frac{d}{dt}  \left\{ \begin{array} {c}
          x   \\
        x_{t1}  \\
        x_{t2}
            \end{array}  \right\} =
 \underbrace{\left\{ \begin{array} {c}
         f^{[a]}(x,x_{t1})  \\
            0 \\
            1
         \end{array}  \right\}}_{\hat{f}^{[1]}} +
 \underbrace{\left\{ \begin{array} {c}
         f^{[b]}(x,x_{t2}) \\
            1 \\
            0
         \end{array}  \right\}}_{\hat{f}^{[2]}}
\end{equation}
with $x_{t1},x_{t2}\in\mathbb{R}$. Note that if the systems
\[
  y^{\prime} = \hat{f}^{[A]}(y),  \qquad
  y^{\prime} = \hat{f}^{[B]}(y)
\]
with $y=(x,x_{t1},x_{t2})$ are solvable, then a splitting method
similar to (\ref{eq:splitting}) can be used, since $x_{t1}$ is constant
when integrating the first equation and $x_{t2}$ is constant when
solving the second one. This, in fact, can be considered as a
generalization of the procedure proposed in \cite{sanzserna96cni} for
time-dependent and separable Hamiltonian systems, and is of
interest if the time-dependent part in $f^{[a]}$ and
$f^{[b]}$ is cheap to compute. Otherwise the overall algorithm may
be computationally costly, since these functions have to be
evaluated $s$ times (the number of stages in (\ref{eq:splitting})) per
time step.

Another disadvantage of this simple procedure is the following.
Suppose that,
when the time is frozen, the function $f$ in (\ref{naut.1}) has
a special structure which allows to apply highly efficient splitting schemes.
If now $t$ is a variable,  with (\ref{naut.2})
this time dependency is eliminated but the structure of the
equation might be modified so that one is bound to resort to more
general and less efficient integrators. This issue has been analyzed in detail in
\cite{blanes06smf}.

A second procedure which avoids the difficulties exhibited by the previous example
consists in approximating the exact solution of
(\ref{naut.1}) or equivalently the flow $\varphi_h$ by the
composition
\begin{equation}   \label{naut.3}
   \psi_{s,h}^{[r]} =
     \varphi_{h}^{[\hat{B}_{s+1}]} \circ
      \varphi_{h}^{[\hat{A}_s]} \circ \varphi_{h}^{[\hat{B}_s]}
       \circ \cdots \circ \varphi_{h}^{[\hat{B}_{2}]} \circ
       \varphi_{h}^{[\hat{A}_1]} \circ
       \varphi_{h}^{[\hat{B}_1]},
\end{equation}
where the maps $\varphi_{h}^{[\hat{A}_i]}$,
$\varphi_{h}^{[\hat{B}_i]}$ are the exact $1$-flows corresponding
to the time-independent differential equations
\begin{equation}   \label{naut.4}
   x^{\prime}  =   \hat{A}_i(x), \qquad\quad
   x^{\prime}  =   \hat{B}_i(x),  \qquad  i=1,2,\ldots
\end{equation}
respectively, with
\begin{equation}  \label{naut5}
      \hat{A}_i(x)   \equiv
      h  \sum_{j=1}^k \rho_{ij} f^{[a]}(x,\tau_j),  \qquad\quad
      \hat{B}_i(x)   \equiv
       h \sum_{j=1}^k \sigma_{ij} f^{[b]}(x,\tau_j).
\end{equation}
Here $\tau_j = t_0 + c_j h$ and the (real) constants
$c_j$, $\rho_{ij}$, $\sigma_{ij}$ are chosen such that $\varphi_h
 = \psi_{s,h}^{[r]} + \mathcal{O}(h^{r+1})$.
 Furthermore, the new schemes, when
 applied to (\ref{naut.1}) with the time frozen, reproduce
 the standard splitting (\ref{eq:splitting}). This is accomplished
 by ensuring that $\sum_{j} \rho_{ij} = a_{i}$ and
 $\sum_{j} \sigma_{ij} = b_{i}$. The $c_j$ coefficients, on the other hand, are typically chosen
 as the nodes of a symmetric quadrature rule of order at least $r$. In particular, if
 a Gauss--Legendre
quadrature rule is adopted, with $k$ evaluations of
$f^{[a]}(x,\tau_j)$ and $f^{[b]}(x,\tau_j)$ a method of order
$r=2k$ can be built (taking $s$ sufficiently large).

Once the
quadrature nodes $\tau_j$ and the number of stages $s$ are fixed,
there still remains to obtain the coefficients $\rho_{ij}$,
$\sigma_{ij}$ such that $\psi_{s,h}^{[r]}$ has the desired order.
This is done by requiring that the composition (\ref{naut.3})
match the solution of (\ref{naut.1}) as given by the Magnus
expansion \cite{blanes08tme}. The task is made easier by noticing that the order
conditions to be satisfied by  $\rho_{ij}$  and $\sigma_{ij}$ are
identical both for linear and nonlinear vector fields. Thus,
the problem for the linear case is solved first and then one generalizes
the treatment to arbitrary nonlinear separable problems.

 The integrators of order four
 and six constructed along these lines in \cite{blanes06smf} are generally
 more efficient than standard splitting methods applied to the
 enlarged system (\ref{naut.2}).

\section{Numerical examples with selected methods}
\label{nexamples}

This section intends to illustrate the relative performance
between different splitting methods, and occasionally we compare
with other standard methods. We consider first a relatively simple
problem where most of the methods previously mentioned can be
used, showing their good features.
 We show the interest of the high order methods when accurate results
 are desired and the improvement which can be achieved when choosing
 a method from the most appropriate family of methods for each problem.
Next, we consider a problem which, due to its very particular
structure, allows to build tailored methods whose performance is
much superior to other splitting methods.

\subsection{The perturbed Kepler problem}

 As a first example, we take the perturbed Kepler problem with
Hamiltonian (\ref{eq.14})
\begin{equation}
H=\frac{\displaystyle1}{\displaystyle2}(p_{1}^{2}+p_{2}^{2})-\frac
{\displaystyle1}{\displaystyle r}-\frac{\displaystyle\varepsilon
}{\displaystyle2r^{3}}\,\left(  1-\alpha\frac{\displaystyle3q_1^{2}%
}{\displaystyle r^{2}}\right),  \label{ne.2}%
\end{equation}
where $r=\sqrt{q_1^{2}+q_2^{2}}$ and the additional parameter $\alpha$ has
been introduced for convenience. This Hamiltonian describes in
first approximation the dynamics of a satellite moving into the
gravitational field produced by a slightly oblate spheric planet.
The motion takes place in a plane containing the symmetry axis of
the planet when $\alpha=1$, whereas $\alpha=0$ corresponds to a
plane perpendicular to that axis \cite{meirovich88moa}.

This simple (but non trivial) example constitutes in fact an excellent test bench
for most of the methods of this paper. Notice that
the system is separable into kinetic
and potential parts, and we can use, for instance, the symmetric
second order method (\ref{leapfrog}) which allows us to get higher
order methods by composition, as given in (\ref{eq:compSS}). On the
other hand, since the system is separable into two solvable parts,
then we can also use methods from Table~\ref{tableAB}, which should show
better performances than methods of the same order considered from
the previous family of methods. In addition,  the kinetic energy is quadratic in momenta, so that
 RKN methods from
Table~\ref{tableRKN} can be used, and one expect a further
improvement. Finally, observe that one may split the system as
\begin{equation}\label{SepNI}
  H = H_{0} + \varepsilon H_I,
\end{equation}
where $H_0$ corresponds to the Kepler problem, which is exactly
solvable.
 The Keplerian part of the Hamiltonian can be solved in
action-angle coordinates, where two changes of variables are
needed. Alternatively, if desired, $H_{0}$ can be integrated in
cartesian coordinates using the $f$ and $g$ Gauss functions, but
then a nonlinear equation must be solved with an iterative
scheme \cite{danby88foc}. In any case, if $\varepsilon\ll 1$, methods from
Table~\ref{tableNI} can be used which should be superior to all
previous methods in the limit $\varepsilon\rightarrow 0$.

We must also mention that the performance of all methods previously mentioned
can be further improved by using the processing technique, and even additional
improvements can be achieved if modified potentials are considered.

We take $\varepsilon=0.001$, which approximately corresponds to a
satellite moving under the influence of the Earth \cite{kirchgraber88aos}
and initial conditions $q_1=1-e$, $q_2=0$, $p_{1}=0$,
$p_{1}=\sqrt{(1+e)/(1-e)}$, with $e=0.2$ (this would be the eccentricity
for the unperturbed Kepler problem). In general, no closed
orbits are present and a precession is observed. Notice that for
the Hamiltonian (\ref{ne.2}) the strength of the perturbation
depends obviously of the value of $\varepsilon$, but also on the
initial conditions. We take $\alpha=1$ and determine numerically
the trajectory for up to the final time $t_{f}=500\cdot2\pi$  (the
exact solution is accurately approximated using a high order
method with a very small time-step, and this computation was
repeated with different time steps and methods to assure the
accuracy is reached up to round off).

 To compare the performance
of different methods it is usual to consider efficiency curves. We
measure the average error in energy computed at times
$t_k=k\cdot2\pi$ for $k=401,402,\ldots,500$ and this is repeated
several times for each method and using different time steps
(changing the computational cost for the numerical integration).

\begin{figure}[h!]
\begin{center}
 \makebox{\epsfig{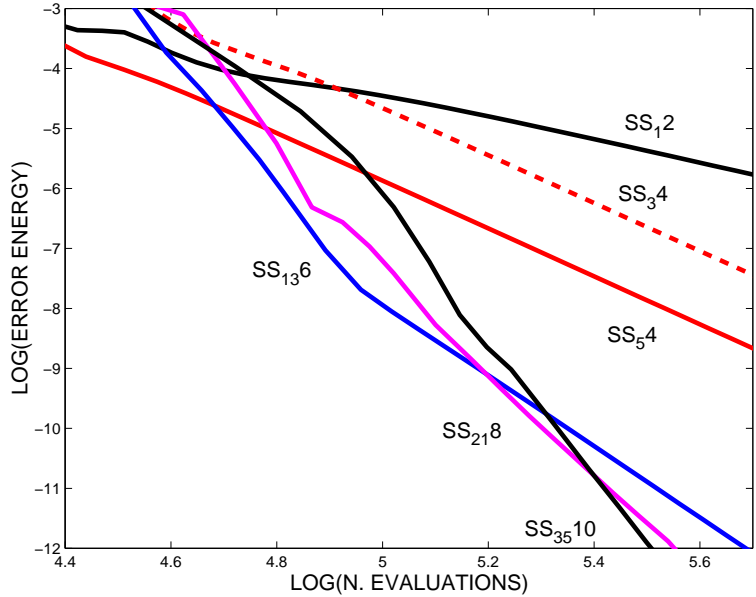}}
\end{center}
\caption{{Average error in energy versus number of force
evaluations in a double logarithmic scale for the numerical
integration of the Hamiltonian system (\ref{ne.2}). It is shown
performance of the most efficient non-processed
symmetric-symmetric methods from Table~\ref{tableSS}.}}
 \label{figEj1-a}
\end{figure}

In the first numerical test, we compare the relative performance
between different symmetric-symmetric methods collected in
Table~\ref{tableSS}.
 We choose as the basic method the symmetric second order composition
 (\ref{leapfrog}) to build higher order methods with the composition
 (\ref{eq:compSS}).
 As mentioned, in general, the performance of the methods of
the same order increase with the number of stages for the methods
in Table~\ref{tableSS}. This is illustrated in  Figure.~\ref{figEj1-a}
 where we show the performance of two 4th-order methods with three
 stages (given by (\ref{suzu1})) and five stages (given by  (\ref{suzukiSS4})).
 The results show that for this problem the five stage method is
 more accurate for all computational costs considered. A similar
 feature is observed for the other methods at higher orders
 (with the exception of the 24-stage 8th-order method which was
 obtained in a different way and the 21-stage methods shows a better performance).
 We choose the best method from Table~\ref{tableSS} at each order
(including the well known three-stages 4th-order method as a
reference) where we
denote by SS$_sr$ the corresponding method of order $r$ using a
$s$-stage composition:

\begin{itemize}
    \item SS$_1$2: The 2nd-order method (\ref{leapfrog}) which has the highest possible stability among
splitting methods.
    \item SS$_3$4 The well known 3-stage 4th-order method (\ref{suzu1}).
    \item SS$_5$4 The 5-stage 4th-order method (\ref{suzukiSS4}) \cite{suzuki90fdo}.
    \item SS$_{13}$6, SS$_{21}$8, SS$_{35}$10: The composition from Table~\ref{tableSS} and
      whose coefficients are given in \cite{sofroniou05dos}.
\end{itemize}

 The results are shown in Fig.~\ref{figEj1-a}, where we
clearly observe that the high order methods have better performance when
high accuracy is desired.

The following numerical experiment intends to illustrate the
interest of the methods designed for problems with some particular
structure. For simplicity, in this numerical test we only consider
fourth-order methods from different families of methods which can
be used on this problem, in order to observe the benefit of tuned
methods for problems with particular structures. The following
methods are considered in addition to SS$_3$4 and SS$_5$4:

\begin{itemize}
    \item S$_6$4: The symmetric 6-stage 4th-order method  for separable
      problems \cite{blanes02psp} from Table~\ref{tableAB}.
    \item RKN$_6$4: The symmetric \textbf{6S}$_{BAB}$ 4th-order method for Nystr\"om
      problems \cite{blanes02psp} from Table~\ref{tableRKN}.
    \item NI(8,4): The 5-stage fourth-order method \textbf{5}(8,4)\textbf{S}$_{BAB}$
      given in \cite{mclachlan95cmi} from  Table~\ref{tableNI}.
    \item RK$_{4}$4: The standard 4-stage 4th-order non-symplectic Runge-Kutta
      methods, used as a reference method.
\end{itemize}

\begin{figure}[h!]
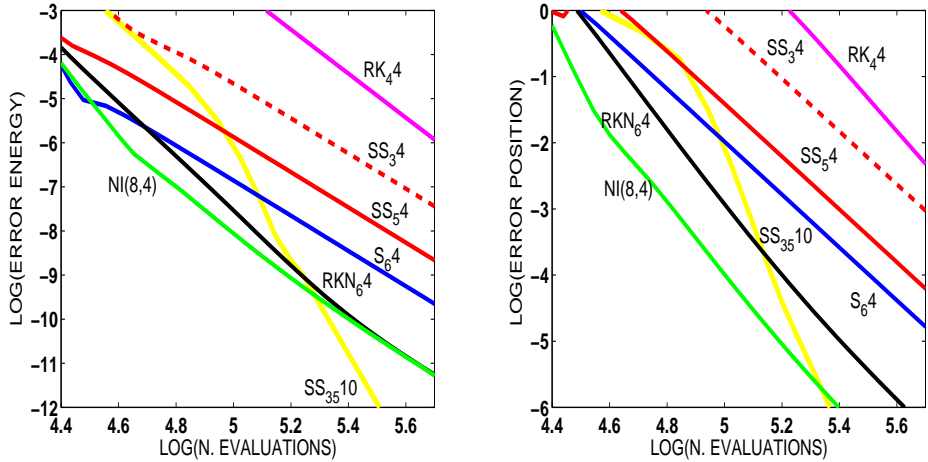

\begin{center}
%\makebox{\epsfig{figure=FigErrEnOrd4.eps,height=6.5cm,width=6.3cm}},
%\makebox{\epsfig{figure=FigErrPosOrd4.eps,height=6.5cm,width=6.3cm}},
\makebox{\epsfig{figure=FigErrEnOrd4b.eps,height=6.5cm,width=6.4cm}}
\makebox{\epsfig{figure=FigErrPosOrd4b.eps,height=6.5cm,width=6.4cm}}
%\makebox{\epsfig{figure=FigErrEnOrd4c.eps,height=6.5cm,width=6.4cm}},
%\makebox{\epsfig{figure=FigErrPosOrd4c.eps,height=6.5cm,width=6.4cm}}
\end{center}
\caption{{Average error in energy (left panel) and position (right
panel) versus number of force evaluations in a double logarithmic
scale for the numerical integration of the Hamiltonian system
(\ref{ne.2}). The performance of different 4th-order methods from
Tables~\ref{tableSS}-\ref{tableNI} is shown. As a reference, we
also shown the results obtained by the standard 4-stage
fourth-order Runge-Kutta method.}}
 \label{figEj1-b}
\end{figure}

Figure~\ref{figEj1-b} shows in double logarithmic scale the
results obtained. In the left panel we show the average error in
energy versus the number of force evaluations and in the right
panel we repeated the same experiment, but we measured the average
error in position (computed at the same instants). For the method
NI(8,4) this counting of the computational cost is not an
appropriate measure. Its computational cost strongly depends on
each particular problem since the evolution of $H_0$ has to be
computed exactly (or very accurately). For simplicity, in our
experiments, we have considered that one stage of NI(8,4) is twice
as expensive as one evaluation of the force. We have also included as a reference
the curve obtained in Fig. \ref{figEj1-a} by SS$_{35}$10.

 Observe that in the first case the results will be
largely independent of $t_{f}$ because the average error in energy
does not increase secularly for symplectic integrators.  For comparison, 
we have also included the results obtained by
the standard 4-stage fourth-order Runge-Kutta method whose error
in energy grows linearly and the error in positions 
quadratically.

Even more accurate results could be obtained as follows. As
mentioned, for this particular problem, modified potentials could
be used and this can be done at a very low computational cost.
Then, methods from Tables~\ref{tableRKNm} and \ref{tableNI} can be
used. For instance, for the split (\ref{SepNI}) we can apply
methods which incorporate modified perturbations $\exp(\varepsilon
C_h(b,c))$ into the algorithm. Then the following map has to be
evaluated:
\begin{align}
\e^{\varepsilon C_h(b,c)}p_{1}  &  =p_{1}+h\varepsilon\left(
b\frac {\displaystyle A}{\displaystyle r^{7}}-h^{2}\varepsilon
c\frac{\displaystyle
C}{\displaystyle(r^{7})^{2}}\right)  q_1 \nonumber\label{ne.3}\\
\e^{\varepsilon C_h(b,c)}p_{2}  &  =p_{2}+h\varepsilon\left(
b\frac {\displaystyle B}{\displaystyle r^{7}}-h^{2}\varepsilon
c\frac{\displaystyle D}{\displaystyle(r^{7})^{2}}\right)  q_2,
\end{align}
where $A=(3/2)(\alpha(3q_1^{2}-2q_2^{2})-r^{2})$,
$B=(3/2)(\alpha5q_1^{2}-r^{2})$,
$C=9(2r^{4}+3\alpha r^{2}(q_2^{2}-4q_1^{2})+\alpha^{2}(18q_1^{4}+q_1^{2}q_2^{2}%
-2q_2^{4}))$ and $D=9(2r^{4}-15\alpha r^{2}q_1^{2}+5\alpha^{2}q_1^{2}(5q_1^{2}%
+2q_2^{2}))$. Notice that the increment in the computational cost
with respect to the evaluation of $\e^{h\varepsilon bF^{[b]}}$
(which corresponds to $c=0$) is only due to a few very simple
additional operations. For this particular example, the evaluation
of the modified perturbation $\e^{\varepsilon C_h(b,c)}$ is about a
$10-20\%$ more expensive than $\e^{hbF^{[b]}}$.

 As a result, more elaborated and
efficient methods are obtained by considering the schemes $(n,4)$
and $(n,5)$ from Table~\ref{tableNI} with processing (which only
requires a few more code lines to program) for RKN problems and
which incorporate modified potentials (see \cite{blanes00psm}).

\subsection{The Schr\"odinger equation}

 As a second example we consider the one-dimensional time-dependent
Schr\"o\-din\-ger equation (\ref{Schrodinger}) with the Morse
potential $V(x)=D \left(1-e^{-\alpha x} \right)^2$. We fix the
parameters to the following values in atomic units (a.u.): $\mu=
1745$ a.u., $ D= 0.2251$ a.u. and $\alpha= 1.1741$ a.u., which are
frequently used for modelling the HF molecule. As initial
conditions we take the Gaussian wave function
%\begin{equation}\label{InGauss}
$  \psi (x,t) = \rho \exp \big( -\beta (x-\bar x)^2 \big)$,
%\end{equation}
with $\beta=\sqrt{k\mu}/2$, $k=2D\alpha^2$, $ \, \bar x=-0.1$ and
$\rho$ is a normalizing constant. Assuming that the system is
defined in the interval $ \, x \in [-0.8, 4.32] $, we split it
into $d=128$ parts of length $ \Delta x=0.04$, take periodic
boundary conditions and integrate along the interval $ \
t\in[0,20\cdot 2 \pi/w_0]$ with $w_0=\alpha \sqrt{2D/\mu}$ (see
\cite{blanes06sso} for more details on the implementation of the
splitting methods to this particular problem).

%  As we have seen, the splitting methods considered in this
%work preserve symplecticity but not unitarity.
% In Figure~\ref{fig4} we show the error in the preservation of
%unitarity, $|q^T(t)q(t)+p^T(t)p(t)-1|$, and the relative error in
%energy ($|(E(t)-E(0))/E(0)|$, where
%$E(t)=u^T(t)Hu(t)=q^T(t)Hq(t)+p^T(t)Hp(t)$) for the 4-stage fourth
%order methods RK$_44$ and GM$_44$. Both methods require the same
%number of FFT calls per step (GM$_44$ has less storage
%requirements) and they are used with the same time step $h=(2
%\pi/w_0)/250$. The error in energy does not grow secularly in time
%for the scheme GM$_44$, as expected from a symplectic integrator,
%whereas the error for the non-symplectic scheme RK$_44$ grows
%linearly. A similar behaviour is observed for the error in
%unitarity, in agreement with the results presented in this work.

%\begin{figure}[h!]
%\begin{center}
%\makebox{\epsfig{figure=fig4mat.eps,height=8cm,width=10cm}}
%\end{center}
%\caption{Error in the preservation of unitarity
% and energy as a function of time in double logarithmic scale
%  for the non symplectic RK$_44$ and symplectic splitting
% GM$_44$ methods applied to the Schr\"odinger equation. Both 4-stage
%  explicit fourth order methods are used with the same time step.}
% \label{fig4}
%\end{figure}

\begin{figure}[h!]
\begin{center}
 \makebox{\epsfig{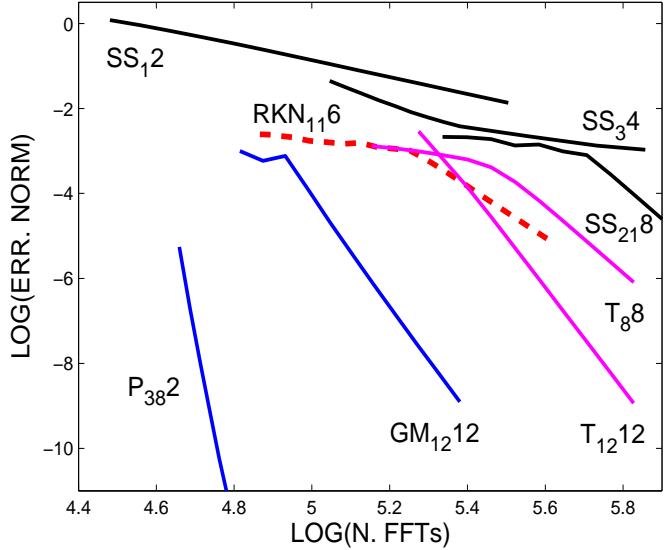}}
\end{center}
\caption{{Error of the vector solution for the Schr\"odinger
equation versus the number of FFT calls in a log-log scale for the
symmetric-symmetric composition methods: SS$_k$n, for methods of
order $n$ using $k$-stage compositions; RKN$_{11}$6 is an 11-stage
6th-order methods from \cite{blanes02psp} for Nystr\"om problems;
the 12-stage 12th-order method, GM$_{12}$12, from
\cite{gray96sit}; and
 P$_{38}2$, a 38-stage second order processed method.}}
 \label{figEj2}
\end{figure}

  Figure~\ref{figEj2}  shows the error in the Euclidean norm of the
vector solution at the end of the integration versus the number of
FFT calls in double logarithmic scale. The integrations are done
starting from a sufficiently small time step and repeating the
computation by slightly increasing the time step until an overflow
occurs, which we identify with the stability limit. We present the
results for the following methods (in addition to the previous
ones SS$_1$2, SS$_3$4 and SS$_{21}$8):
\begin{itemize}
    \item RKN$_{11}$6: The 11-stage 6th-order method \textbf{11S}$_{BAB}$-\cite{blanes02psp}
      from Table~\ref{tableRKN}.
    \item GM$_{12}$12: The 12-stage 12th-order method from \cite{gray96sit}
    tailored  for linear problems with this particular structure.
    \item P$_{38}2$: The 38-stage second order processed method
      with coefficients given in \cite{blanes06sso} tailored  for linear problems with this
      particular structure (only the
computational cost required to evaluate the kernel has been taken
into consideration).
    \item T$_rr, \ r=8,12$: $r$-stage $r$th-order Taylor methods
      obtained by truncating the exponential up to order $r$.
\end{itemize}

From the figure we observe that, for this numerical experiment,
standard Taylor methods outperform to general symmetric-symmetric
splitting methods and are also more accurate than the RKN method.
This is because this problem has a very particular structure and
these splitting methods are not optimized for them. However, the
schemes GM$_{12}$12 and P$_{38}2$ are built for linear problems
with this structure and their superiority is clearly manifest. It
is important to remember that Taylor methods are non-geometric
integrators. The numerical experiments are carried out for a
relatively short time, and the relative performances of the Taylor
methods deteriorates with respect to splitting methods for longer
integrations.

\section{Conclusions and outlook}

Splitting methods are a  flexible and powerful way to solve numerically the
initial value problem defined by (\ref{eq.1.1}) when $f$ can be decomposed into two or more parts and
each of them is simpler to integrate than the original problem.
This is especially true when the exact flow possesses some
structural features which seems natural to reproduce at the discrete level, as happens, for instance, in
Hamiltonian, Poisson, volume-preserving or time-reversible dynamical systems. They are
explicit, usually simple to apply in practice and constitute an important class of geometric numerical
integrators. Closely connected with splitting schemes are composition methods. In this case, the idea
is to construct numerical integrators of arbitrarily high order by composing one or more basic schemes of low order
with appropriately chosen coefficients. The resulting method inherits the relevant properties that the
basic integrator shares with the exact solution, provided these properties are preserved by
composition.

In this paper we have reviewed some of the main features of splitting and composition methods in the
numerical integration of ordinary differential equations. We have presented a novel approach to get the order
conditions of this class of schemes based on Lyndon words and we have seen how these order conditions
particularize when coping with special classes of dynamical systems (near-integrable systems and
second-order differential equations of the form $y''=g(y)$). It turns out that the number of equations to
be solved increases dramatically with the order considered, as so does the complexity of the problem of
finding efficient high order methods. One way to circumvent (up to a certain point) this difficulty consists in
applying the processing technique, since then it is possible to design algorithms with fewer evaluations per
time step. In this sense, one could say that the use of processing is perhaps the most economical path
to achieve high order.

Since splitting methods are widely applied in many areas of science, it is not surprising that a great number
of different schemes are available in the mathematical, physical and chemical literature. We have collected
here some of the most representative integrators, classified according to the particular structure of the differential
equations, the number of stages and the order of consistency, citing in each case the actual reference
where the method has been first proposed.

The good qualitative behavior exhibited by splitting methods  (including preservation
of invariants and structures in phase space), as well as their favorable error propagation in long-time integrations
can be accounted for by applying the theory of backward error analysis. Loosely speaking, the observed
performance is related with the fact that the numerical solution provided by the splitting method is the exact
solution of a differential equation with the same geometric properties as the original system. This
interpretation constitutes in addition the basis for rigorous estimates on the numerical solution.

In contrast with standard integration methods (Runge--Kutta, multistep), whose efficiency is essentially
independent of the particular differential equation considered, splitting schemes can be designed
to incorporate in their formulation some of the most relevant properties of the original system. This
feature has to be taken into account when comparing the efficiency of splitting methods with respect to other
general purpose integrators. In this sense, Figures \ref{figEj1-a} and \ref{figEj1-b} are quite illustrative.
For this particular problem,
specially adapted 4th-order splitting schemes are  up to 6 orders of magnitude
more accurate with the same computational cost than the well known Runge--Kutta method. They even
outperform other standard higher order composition integrators for a wide range of values of the step size
$h$.

As an additional evidence of the extraordinary flexibility of splitting methods, we have considered
the problem of designing specially tailored schemes for the numerical integration of the generalized
harmonic oscillator (\ref{harmonicGeneral}). It turns out that several partial differential equations
appearing in quantum mechanics, optics and electrodynamics give rise, once discretized in space, to
this system with different matrices $M$ and $N$. The particular structure of this dynamical system
can be exploited to build an optimized
 processed second order method involving a large number of stages that
nevertheless is far more efficient than other integrators.

There are other issues in connection with splitting and composition methods that
we have \emph{not} tackled here,
however, and that are also important in this context. Among them we can mention the following.
\begin{itemize}
 \item As was remarked in the introduction, no general rule is provided here to split
  any given function $f$ in the differential equation (\ref{eq.1.1}). It turns out that, for
  $f$ within a certain class of ODEs, this can be done systematically, whereas for other functions
  one has to proceed on a case by case basis. Sometimes, several splittings are possible, and
  the different schemes built from them lead to the preservation of distinctive geometric
  properties. It makes sense, then, to classify the ODEs and their corresponding integration
  methods into different categories. This aspect has been analyzed in \cite{mclachlan02sm}.
  Moreover, in many physical problems there are several geometric properties that are
  conserved simultaneously along the evolution and it is not clear at all how to design
  methods preserving all of them. In that case, which one is the most relevant from a numerical
  point of view?
 \item In this paper we have only considered the initial value problem defined by eq.
 (\ref{eq.1.1}) and integration methods with constant step size $h$. Backward error analysis
 provides an argument why this has to be the case in geometric numerical integration: the
 modified equation corresponding to the numerical method depends explicitly on $h$, so that
 if $h$ is changed so does the modified equation and no preservation of geometric properties
 is guaranteed. There are problems, however, where the use of an adaptive step size is
 mandatory, for instance in configurations of the $N$-body problem allowing close encounters.
 In this case one may apply splitting methods with variable step size by using some specifically
 designed transformations involving the time variable, in such a way that in the new variables
 the resulting time step is constant (see, e.g., \cite{blanes05agi}).
\item  As we have shown in section \ref{ntsteps}, the presence of negative coefficients in
splitting methods of order higher than two is unavoidable. This is not a problem when the flow
of the differential equation evolves in a group (such as in the Hamiltonian case), but may be
unacceptable when the ODE originates from a partial differential equation that is ill-posed
for negative time progression. Several alternatives have been proposed in the literature,
mainly by considering, when possible, modified potentials \cite{blanes05otn}, as noted in
section \ref{modified}. One should
observe, however, that the analysis done in section \ref{ntsteps} does not preclude the
existence of \emph{complex} coefficients with positive real parts. As a matter of fact, splitting methods
with complex coefficients have been developed and tested for problems
in which the Hamiltonian is split into kinetic and potential energy terms \cite{chambers03siw}, for the
time-dependent Schr\"odinger equation \cite{bandrauk06cis}, for generic parabolic equations
\cite{castella08smw} and also in the more abstract
setting of evolution PDEs in analytic semigroups \cite{hansen08hos}.
\item An important characteristic of any numerical integration method  is
  \emph{stability}. Roughly speaking, the
numerical solution provided by a stable numerical
integrator does not tend to infinity when the exact solution is bounded.
Although important, this feature has received considerably less attention in the
specific case of splitting methods. To test the (linear) stability of the method (\ref{eq:splitting}),
instead of the linear equation $y^{\prime} = a y$ as in the usual
stability analysis for ODE integrators, one considers the harmonic oscillator
$y^{\prime\prime} + \lambda^2 y = 0$, $\lambda >0$, as a model problem with a splitting of the
form (\ref{harmonic2}). The
idea is to find the time steps for which all numerical solutions
remain bounded.  The integrator (\ref{eq:splitting}) typically will be
unstable for $|h \lambda| > x_*$, where the parameter $x_*$
determines the stability threshold of the numerical scheme. In particular, for the
leapfrog method one has $x_* = 2$. Although the stability threshold imposes restrictions on the
step size, in the process of building
high order schemes, linear stability  is not usually taken into account, ending
sometimes with methods possessing such a small relative stability threshold that they are
useless in practice. In this way, constructing high order splitting methods with relatively large
linear stability intervals and highly accurate is of great interest. This has been achieved in
reference \cite{blanes07otl} for linear systems, but remains an open problem in general.

\item In section \ref{sec.5} we have mentioned an optimization criterion to choose
the free parameters in splitting and composition methods, which consist in minimizing
the Euclidean norm of the coefficients that constitute the leading error term of the method. It is
clear, however, that minimizing the leading error term does not guarantee that the method thus
obtained is the most efficient: it might occur that the influence of the subsequent error terms
is the decisive factor in the performance of the scheme. In this sense, it would be extremely interesting
to have estimates on all the error terms in the asymptotic expansion of the modified equation and get
the coefficients of the method that minimize these estimates.

\item The numerical analysis of second-order differential equations with oscillatory
solutions has aroused much interest during the past few years.  The typical test problem in this
setting is the equation $q'' + \Omega^2 q = f(q)$, where $\Omega$ is a symmetric and positive definite
matrix. Here the aim is to design new methods which improve in accuracy and stability the standard
St\"ormer--Verlet integrator. We refer the reader to \cite{hairer06gni} and references therein for a comprehensive study of this problem.

\item Although only ODEs have been considered here, splitting methods have been
also applied with success to stochastic differential equations (SDEs). Here the
aim is, as in the deterministic case, to design integration
methods which automatically incorporate conservation properties the SDE  possesses
\cite{misawa01ala}.

\end{itemize}

\subsection*{Acknowledgements}

This work has been partially supported by Ministerio de
Ciencia e Innovaci\'on (Spain) under project MTM2007-61572
(co-financed by the ERDF of the European Union). SB also aknowledges 
the support of the UPV through the project 20070307.
We are especially grateful to Prof. Arieh
Iserles for his kind invitation to write and submit this paper.

\bibliographystyle{plain}
\bibliography{ourbib,geom_int,numerbib}

\end{document}